\documentclass[10pt]{article}
\usepackage[textwidth=468pt,textheight=625pt,centering]{geometry}
\usepackage[numbers,sort&compress]{natbib}
\usepackage{amssymb}
\usepackage{amsmath}
\usepackage{bm}
\usepackage{multirow}
\usepackage{amsthm}
\usepackage{caption}
\usepackage{subcaption}
\usepackage{array}
\usepackage{xcolor}
\usepackage{graphicx}
\usepackage{xltabular}

\def\R{\mathbb{R}}

\def\ie{i.e.}
\def\eg{e.g.}

\newcommand{\bracket}[1]{\left( #1 \right)}
\newcommand{\set}[1]{\left\{ #1 \right\}}
\newcommand{\defset}[2]{\left\{ #1: \ #2 \right\}}
\newcommand{\abs}[1]{\left| #1 \right|}

\newcommand{\tnorm}[1]{\|#1\|}
\newcommand{\C}[1]{{C}^{#1}}
\newcommand{\poly}{\mathbb P}
\newcommand{\PSR}[1]{{#1}_{\mathrm{PS}}}
\newcommand{\PST}{\PSR{\triangle}}
\newcommand{\vertices}{V}
\newcommand{\bvertices}{\vertices^b}
\newcommand{\edges}{E}
\newcommand{\bedges}{\edges^b}
\newcommand{\triangles}{T}
\newcommand{\conv}[1]{[ #1 ]}
\newcommand{\convhalf}[1]{[ #1 )}
\newcommand{\grad}[1]{\nabla #1}

\newtheorem{example}{Example}

\begin{document}

\title{Isogeometric analysis with $C^1$ cubic Powell--Sabin splines}

\author{Jan Gro\v{s}elj$^{1,2}$ \and Ada \v{S}adl Praprotnik$^{1,2,*}$ \and Hendrik Speleers$^{3}$}

\date{}

\maketitle

{\small
\noindent $^{1}$ Faculty of Mathematics and Physics, University of Ljubljana, Jadranska 19, 1000 Ljubljana, Slovenia\\
\noindent $^{2}$ Institute of Mathematics, Physics and Mechanics, Jadranska 19, 1000 Ljubljana, Slovenia\\
\noindent $^{3}$ Department of Mathematics, University of Rome Tor Vergata,
Via della Ricerca Scientifica 1, 00133 Rome, Italy\\
\noindent $^{*}$ Corresponding author. E-mail: ada.sadl-praprotnik@fmf.uni-lj.si
}

\begin{abstract}
In this paper, we consider $\C{1}$ cubic Powell--Sabin splines for the numerical solution of boundary value problems on planar and spatial surface domains. We first review the construction and basic properties of polynomial and rational $\C{1}$ cubic Powell--Sabin spline representations on unstructured triangulations. Then, we discuss how these flexible representations can be exploited to create geometry mappings suited for a precise description of (classes of) surface domains. This is illustrated with several examples. Finally, the obtained domain descriptions are utilized in the isogeometric analysis framework for solving various Poisson and biharmonic problems. It is demonstrated that $\C{1}$ cubic Powell--Sabin splines form a powerful alternative to $\C{0}$ cubic Lagrange elements and bicubic NURBS.
\end{abstract}

domain parametrization

\section{Introduction}
\label{sec:intro}

Polynomial splines over a triangulation are an extensively studied approximation framework in numerical analysis. Numerous constructions have their origin in the theory of finite elements, where they serve as building blocks in space discretizations for solving partial differential equations (PDEs). In this context, they are traditionally characterized by an interpolation problem such that each polynomial piece over a triangle of the triangulation is determined by local interpolation data \cite{fem_ciarlet_02}.

Isogeometric analysis is a spin-off of finite element analysis focusing mainly on bridging the gap between simulation and design \cite{iga_hughes_05,iga_cottrell_09}. In design systems, models are commonly parametrized by tensor-product B-splines, \ie{}, a set of locally supported basis functions that form a convex partition of unity \cite{deboor_01}. The splines expressed in such a basis come equipped with a geometrically intuitive control structure and stable evaluation algorithms. The key concept of isogeometric analysis is to use these basis functions also as an approximation framework to solve PDEs in the isoparametric sense. This means that both, the parametrization of the domain and the numerical solution on the domain, are expressed by the same set of basis functions. For this approach to be computationally effective, a crucial step has been the generalization of standard basis functions to support adaptive refinement strategies; see, \eg{}, \cite{ts_sederberg_03,thb_giannelli_12,lr_dokken_13}.

Although finite element analysis and isogeometric analysis are built on the same foundations, the theories have quite ambivalent outlook on the underlying mesh and space discretizations. Isogeometric analysis highlights the role of splines in the approximation process and is focused on the investigation of effects that spline properties, such as high-order smoothness and convex representation, have on the accuracy and stability of the numerical solution of a PDE. On the other hand, finite element analysis is commonly performed on geometries discretized by triangular meshes using only continuous or even discontinuous splines that are considered separately on each element of the partition.

In order to adapt splines over a triangulation to the isogeometric analysis framework, an important step is to establish spline representations in terms of basis functions that mimic the properties of B-splines, in particular the nonnegativity and partition of unity. In a local sense, such a representation is provided by the Bernstein--B\'{e}zier \cite{farin_86,lai_07} or simplex spline \cite{cohen_13,lyche_18,lyche_22,lyche_24,lyche_25} techniques. The transition between local and global context can be carried out by finding a minimal determining set of the chosen smoothness, which identifies the degrees of freedom and presents an alternative to using interpolation functionals. This approach has been successfully tested in applications related to isogeometric analysis; see, \eg{}, \cite{iga_jaxon_14,iga_xia_15,iga_wang_18,iga_xia_18,iga_zareh_19,iga_eddargani_24}.

Global B-spline techniques over a triangulation are looking for locally supported basis functions that form a convex partition of unity and have a certain order of smoothness on the entire domain. Such basis functions can be constructed by similar techniques as local basis functions but usually require specific configurations or geometric assumptions that allow for the direct incorporation of smoothness properties into their definition. An example of such basis functions obtained by means of simplex splines are the so-called TCB-splines \cite{liu_07,wang_22}. On the other hand, the Bernstein--B\'{e}zier techniques have enabled the development of several B-spline-like representations for classical and non-classical finite elements of Argyris \cite{bell_groselj_20,arg_groselj_21}, Clough--Tocher \cite{rct3_speleers_10,ct_groselj_22}, and Powell--Sabin \cite{ps_dierckx_97,ps5_speleers_10,ps_speleers_13,ps_groselj_16} type.

This paper is focused on utilizing $\C{1}$ cubic Powell--Sabin B-splines in the isogeometric analysis framework. The B-spline representation for Powell--Sabin splines of degree~$2$ is well-established and widely tested in the context of geometric modeling \cite{ps_windmolders_99,ps_windmolders_00,rps_speleers_07,rps_speleers_13} and finite element applications \cite{speleers_06,speleers_07,speleers_12,ps2_speleers_15,iga_giorgiani_18}. The degree~$3$ offers several options for definition of a suitable spline space \cite{ps3_chen_08,ps3_lamnii_14,ps3_speleers_15,ps3_groselj_16,ps3_groselj_17}. Here, we consider the macro-element space introduced in \cite{ps3_speleers_15} that has a dimension very similar to $\C{0}$ cubic Lagrange elements, possesses some attractive $\C{2}$ super-smoothness properties, and is globally representable by a set of B-spline-like functions. We take advantage of the properties of these functions and generalize the piecewise polynomial space to a vector space spanned by rational B-spline functions.

A rational B-spline representation over an unstructured triangulation provides a versatile approximation framework. In the context of isogeometric analysis it can be used to express, either in exact or approximate way, the geometries that are usually described by tensor-product NURBS. This may gracefully resolve some common issues with NURBS geometries, such as $\C{1}$ spline parametrization at extraordinary points, irregular local mesh refinements, and restriction to four-sided domains.

In this paper, we develop techniques and present examples of geometry descriptions based on $\C{1}$ cubic Powell--Sabin B-splines. By studying the properties of the basis functions, we also discuss how to satisfy boundary conditions imposed on a domain. Then, we apply the splines to solve second and fourth order boundary value problems on surface domains with the isogeometric analysis approach. Through numerical examples we demonstrate that these splines perform comparably to or better than $\C{0}$ cubic Lagrange elements and bicubic NURBS, thereby providing a powerful tool combining the benefits of the two standard frameworks.

The remainder of the paper is organized as follows. Section~\ref{sec:ps3_splines} provides an overview of the $\C{1}$ cubic Powell--Sabin B-splines and introduces their rational generalization. In Section~\ref{sec:approx}, techniques for representing and constructing geometry mappings are discussed. This section also provides a setup for solving boundary value problems and the treatment of boundary conditions. Section~\ref{sec:problems} contains several numerical examples of solving Poisson and biharmonic problems on surface domains. In Section~\ref{sec:conclusion}, some concluding remarks are presented.

\section{Cubic Powell--Sabin splines}
\label{sec:ps3_splines}

\subsection{Powell--Sabin refinement}
\label{sec:ps_refinement}

Let $\Theta \subset \mathbb{R}^2$ be a closed polygonal region with connected boundary $\partial \Theta$. We consider this region as a parametric domain. Moreover, let $\triangle$ be a triangulation of $\Theta$ represented by a triplet $\triangle = (\vertices, \edges, \triangles)$, where $\vertices$ is the set of vertices, $\edges$ is the set of edges, and $T$ is the set of triangles of the triangulation $\triangle$. These sets are more precisely specified as follows.
\begin{itemize}
\item The set $\vertices$ consists of points $v_i \in \Theta$ indexed by $i$ in range from $1$ to the cardinality $|\vertices|$ of $\vertices$.

\item The set $E$ consists of edges $e_l$ indexed by $l$ in range from $1$ to the cardinality $|\edges|$ of $\edges$. An edge is the convex hull $\conv{v_i, v_j} \subset \Theta$ of two distinct points $v_i, v_j \in \vertices$.

\item The set $T$ consists of triangles $t_m$ indexed by $m$ in range from $1$ to the cardinality $|\triangles|$ of $\triangles$. A triangle is the convex hull $\conv{v_i, v_j, v_k} \subseteq \Theta$ of three noncollinear points $v_i, v_j, v_k \in \vertices$.
\end{itemize}
The union of triangles in $\triangles$ is $\Theta$, and we impose that the intersection of two distinct triangles from $\triangles$ is either empty or corresponds to a common vertex or edge of the two triangles. We also assume that each edge in $\edges$ is an edge of a triangle in $\triangles$, and each vertex in $\vertices$ is a vertex of a triangle in $\triangles$. Additionally, we denote by $\bvertices \subseteq \vertices$ the set of boundary vertices and by $\bedges \subseteq \edges$ the set of boundary edges. Namely, $v_i \in \vertices$ is in $\bvertices$ if $v_i \in \partial \Theta$, and $e_l \in \edges$ is in $\bedges$ if $e_l \subset \partial \Theta$.

The Powell--Sabin refinement of $\triangle$ is a triangulation $\PST$ obtained by splitting each triangle of $\triangles$ into six smaller triangles in the following way. Let $t_m \in \triangles$. Suppose $e_l \in \edges \setminus \bedges$, and let $t_{m'} \in \triangles$ be the triangle sharing the edge $e_l$ with $t_m$. Choose the triangle split points $v_m^t$ and $v_{m'}^t$ in the interior of $t_m$ and $t_{m'}$, respectively, so that the line segment $\conv{v_m^t, v_{m'}^t}$ intersects the interior of the edge $e_l$. The point of intersection is denoted by $v_l^e$ and is called an edge split point. If $e_l \in \bedges$, choose the split point $v_l^e$ to be any point in the interior of the edge. The refinement of $t_m$ into six smaller triangles is obtained by connecting the triangle split point $v_m^t$ to the vertices of $t_m$ and to the edge split points located on the edges of $t_m$.

\begin{figure}[!t]
\centering
\includegraphics[width=0.6\textwidth]{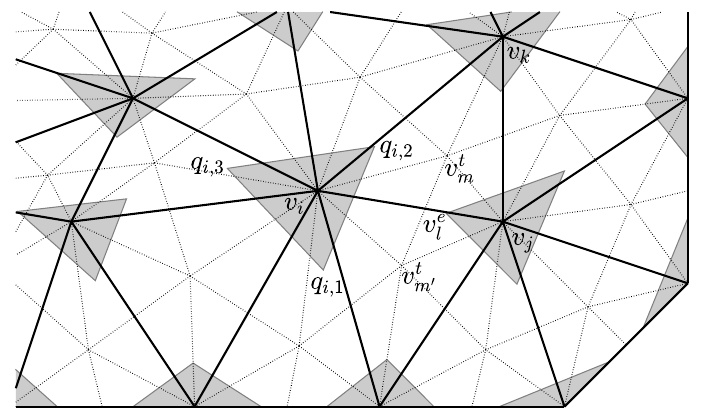}
\caption{An excerpt from a triangulation $\triangle$ (depicted by solid lines) that is refined by a Powell--Sabin triangulation $\PST$ (depicted by dotted lines). The gray colored triangles corresponding to the vertices of $\triangle$ represent a configuration for construction of cubic Powell--Sabin B-splines.}
\label{fig:pstri}
\end{figure}

Applying the refinement process to each triangle of $\triangles$ produces the Powell--Sabin refinement $\PST$ of the original triangulation $\triangle$. An example of the refinement is shown in Figure~\ref{fig:pstri}. Note that the geometric condition involved in the choice of the triangle split points can be automatically satisfied by choosing the split points to be the incenters of the corresponding triangles. We represent the refined triangulation by $\PST = (\PSR{\vertices}, \PSR{\edges}, \PSR{\triangles})$, where $\PSR{\vertices}$ is the set of vertices, $\PSR{\edges}$ is the set of edges, and $\PSR{\triangles}$ is the set of triangles of $\PST$. Each $\PSR{t} \in \PSR{\triangles}$ is of the form $\PSR{t} = \conv{v_i, v_l^e, v_m^t}$ for some indices $i$, $l$, $m$, such that $t_m \in \triangles$, $e_l \in \edges$ is an edge of $t_m$, and $v_i \in \vertices$ is an endpoint of $e_l$.

\subsection{Cubic B-spline representation}
\label{sec:PS_Bspline}

Let $\poly_3$ denote the space of bivariate polynomials of total degree at most three. For a given triangulation $\triangle$ of $\Theta$ and a refinement $\PST$ of $\triangle$ we consider the vector space
\begin{equation*}
\mathbb{S}_3(\PST) =
\defset{s \in \C{1}(\Theta)}{s|_{\PSR{t}} \in \poly_3, \ \PSR{t} \in \PST} \cap
\C{2}\bracket{\defset{v_m^t}{t_m \in \triangles}} \cap
\C{2}\bracket{\defset{v_l^e}{e_l \in \bedges}},
\end{equation*}
\ie{}, the space of cubic splines over $\PST$ that are $\C{1}$ smooth everywhere and $\C{2}$ smooth at every triangle split point and at every boundary edge split point. Remarkably, these additional $\C{2}$ super-smoothness conditions imply that every spline in $\mathbb{S}_3(\PST)$ is $\C{2}$ smooth at $\convhalf{v_m^t, v_l^e}$ for any triangle $t_m \in \triangles$ and any edge $e_l \in \edges$ attached to this triangle.

\begin{figure}[!t]
\centering

\begin{subfigure}[b]{\textwidth}
\begin{subfigure}[b]{0.32\textwidth}
\centering
\includegraphics[width=\textwidth]{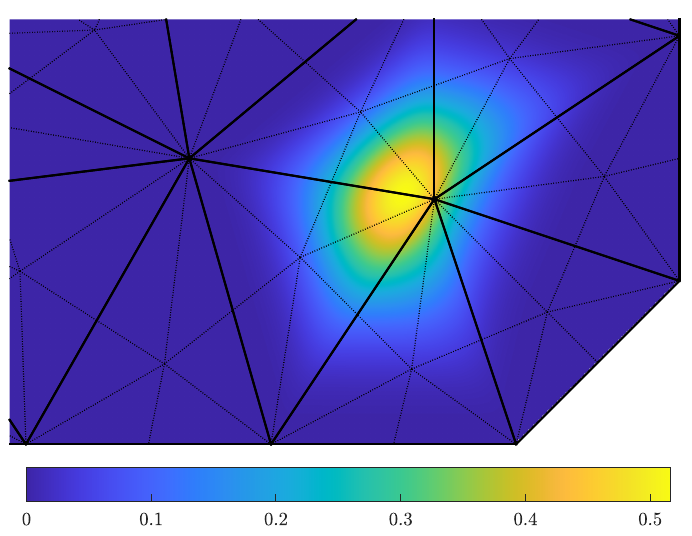}
\end{subfigure}
\hfill
\begin{subfigure}[b]{0.32\textwidth}
\centering
\includegraphics[width=\textwidth]{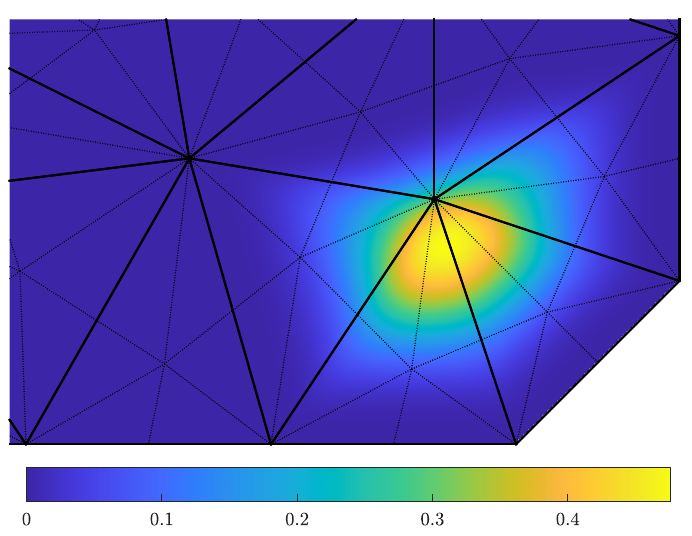}
\end{subfigure}
\hfill
\begin{subfigure}[b]{0.32\textwidth}
\centering
\includegraphics[width=\textwidth]{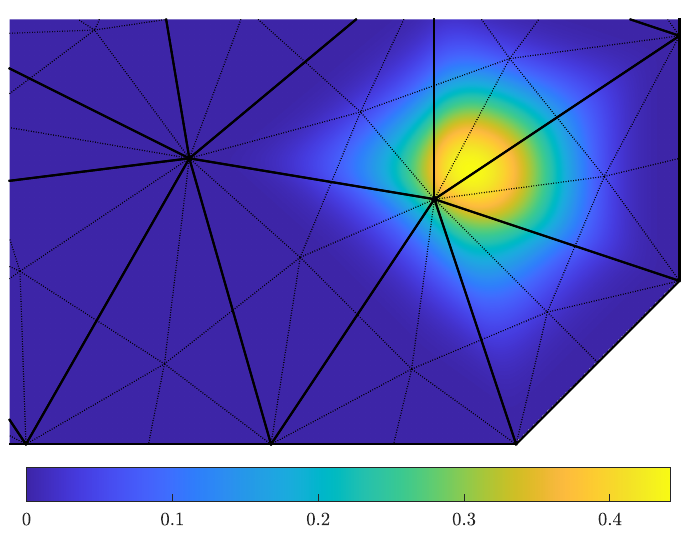}
\end{subfigure}
\caption{Three B-splines associated with an interior vertex.}
\label{fig:psbv1}
\end{subfigure}

\begin{subfigure}[b]{\textwidth}
\begin{subfigure}[b]{0.32\textwidth}
\centering
\includegraphics[width=\textwidth]{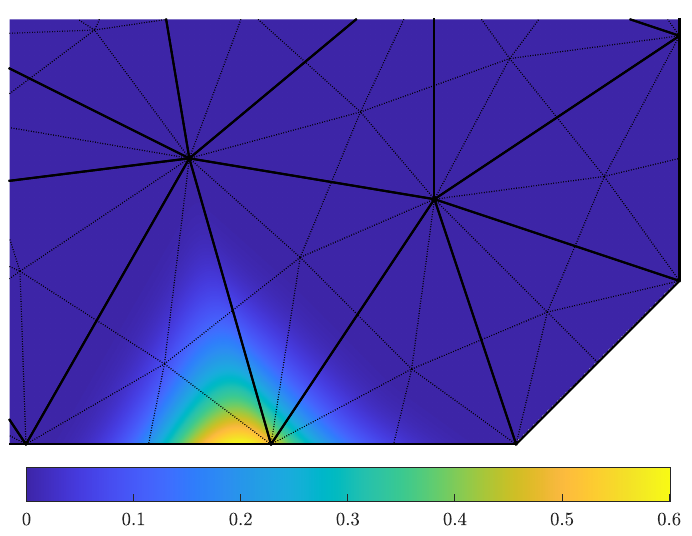}
\end{subfigure}
\hfill
\begin{subfigure}[b]{0.32\textwidth}
\centering
\includegraphics[width=\textwidth]{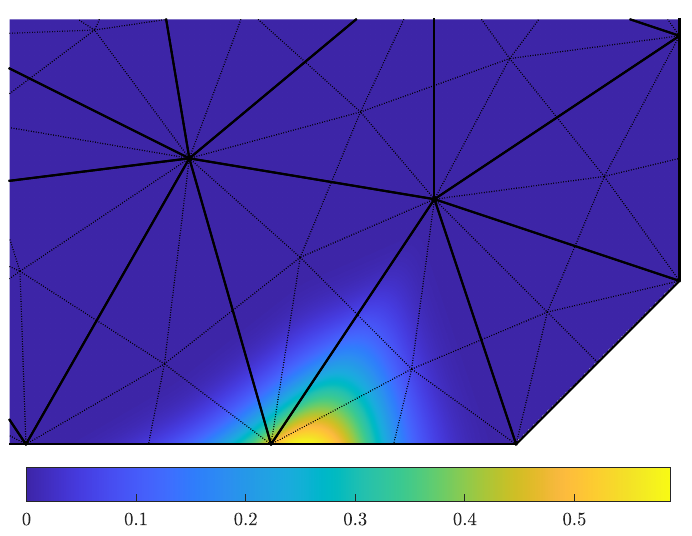}
\end{subfigure}
\hfill
\begin{subfigure}[b]{0.32\textwidth}
\centering
\includegraphics[width=\textwidth]{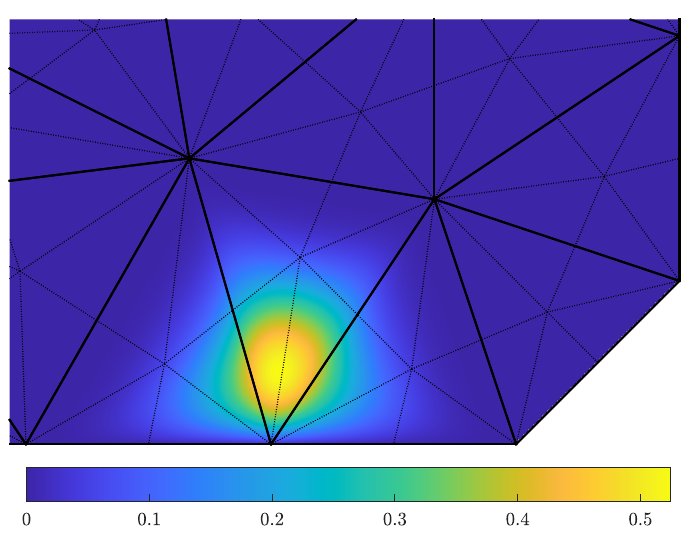}
\end{subfigure}
\caption{Three B-splines associated with a boundary vertex.}
\label{fig:psbv2}
\end{subfigure}

\begin{subfigure}[b]{\textwidth}
\begin{subfigure}[b]{0.32\textwidth}
\centering
\includegraphics[width=\textwidth]{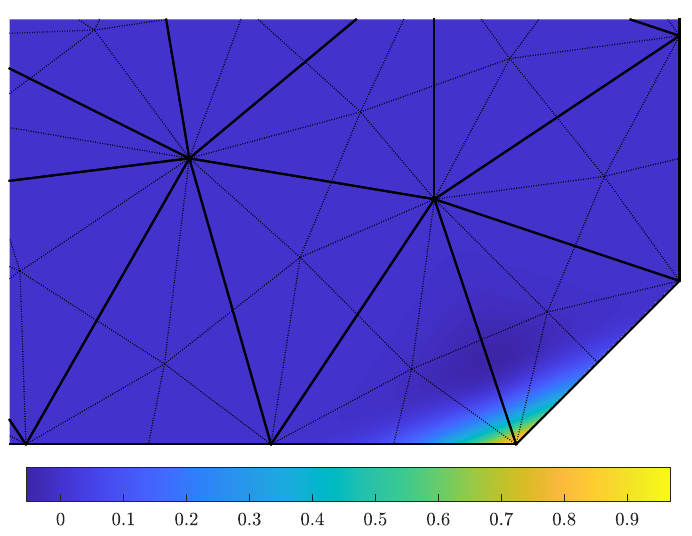}
\end{subfigure}
\hfill
\begin{subfigure}[b]{0.32\textwidth}
\centering
\includegraphics[width=\textwidth]{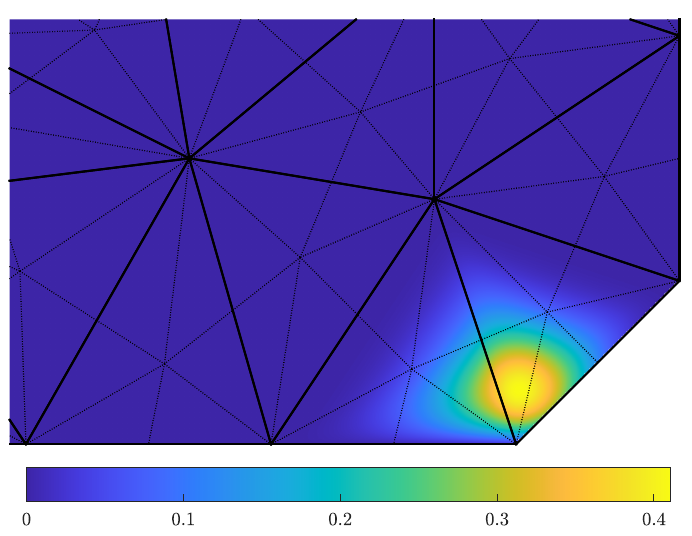}
\end{subfigure}
\hfill
\begin{subfigure}[b]{0.32\textwidth}
\centering
\includegraphics[width=\textwidth]{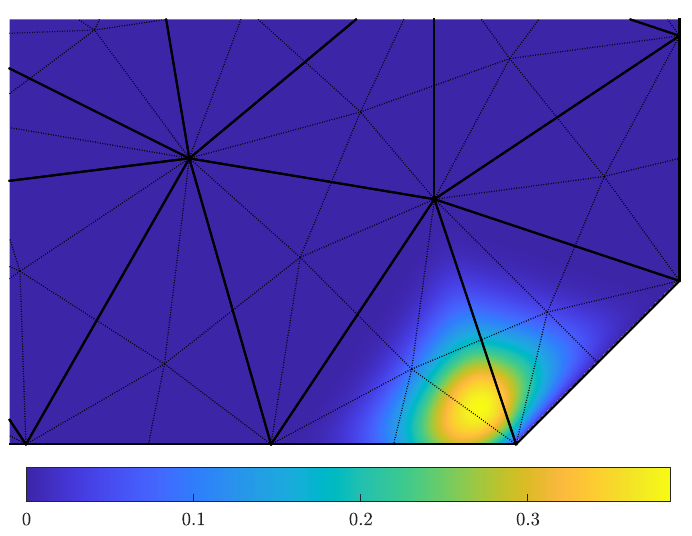}
\end{subfigure}
\caption{Three B-splines associated with a corner vertex.}
\label{fig:psbv3}
\end{subfigure}

\caption{Examples of cubic Powell--Sabin B-splines associated with vertices of the triangulation shown in Figure~\ref{fig:pstri}.}
\label{fig:psbv}
\end{figure}

\begin{figure}[!t]
\centering

\begin{subfigure}[b]{0.49\textwidth}
\begin{subfigure}[b]{0.48\textwidth}
\centering
\includegraphics[width=\textwidth]{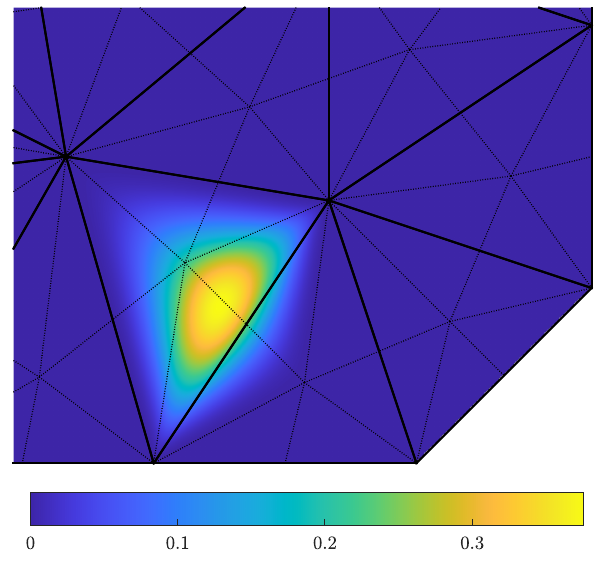}
\end{subfigure}
\hfill
\begin{subfigure}[b]{0.48\textwidth}
\centering
\includegraphics[width=\textwidth]{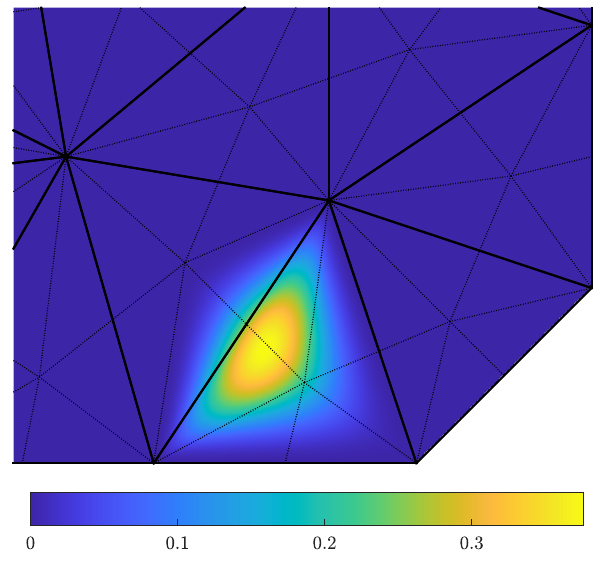}
\end{subfigure}
\caption{Two B-splines associated with an interior edge.}
\label{fig:psbe1}
\end{subfigure}
\hfill
\begin{subfigure}[b]{0.49\textwidth}
\begin{subfigure}[b]{0.48\textwidth}
\centering
\includegraphics[width=\textwidth]{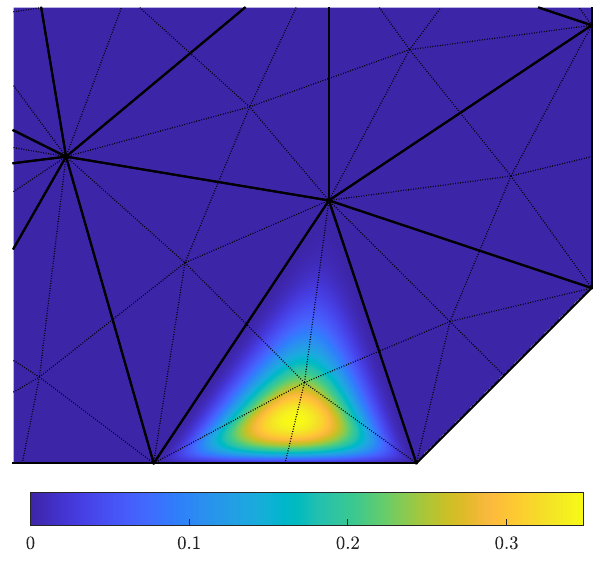}
\end{subfigure}
\hfill
\begin{subfigure}[b]{0.48\textwidth}
\centering
\includegraphics[width=\textwidth]{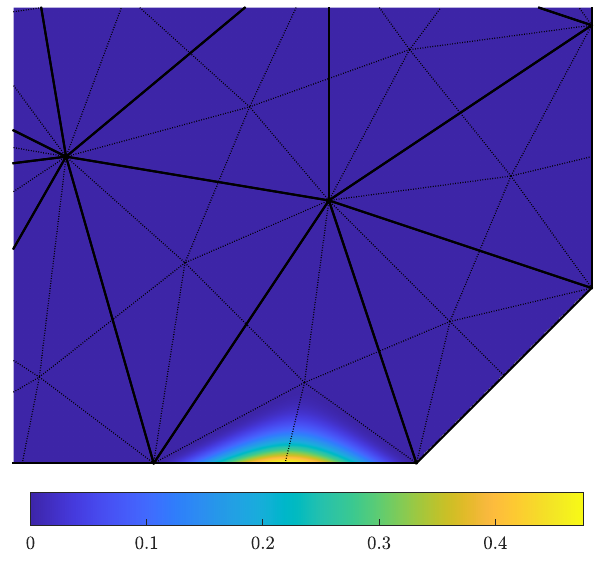}
\end{subfigure}
\caption{Two B-splines associated with a boundary edge.}
\label{fig:psbe2}
\end{subfigure}

\caption{Examples of cubic Powell--Sabin B-splines associated with edges of the triangulation shown in Figure~\ref{fig:pstri}.}
\label{fig:psbe}
\end{figure}

The space $\mathbb{S}_3(\PST)$ was introduced in \cite{ps3_speleers_15} and has a convenient dimension formula, \ie{}, $\dim(\mathbb{S}_3(\PST)) = 3 |\vertices| + 2 |\edges|$. By definition, it contains $\poly_3$ and thereby allows optimal approximation order. Moreover, a spline from $\mathbb{S}_3(\PST)$ can be locally characterized and represented by B-spline-like basis functions. Namely,
\begin{itemize}
\item three basis functions $B_{i,r}^v: \Theta \rightarrow \R$, $r = 1, 2, 3$, for each $v_i \in \vertices$,
\item two basis functions $B_{l,r}^e: \Theta \rightarrow \R$, $r = 1, 2$, for each $e_l \in \edges$.
\end{itemize}
All basis functions are nonnegative and form a partition of unity. The nonnegativity and shape properties are closely related to a choice of triangles associated with the vertices of $\triangle$, which define a configuration of the basis functions (for an illustration, see Figure~\ref{fig:pstri}). This configuration is discussed in Section~\ref{sec:approx_boundary}. Further details on the construction of the basis and its properties are presented in \cite{ps3_speleers_15,ps3_groselj_17}. Examples of these basis functions are shown in Figures \ref{fig:psbv} and \ref{fig:psbe}.

Thanks to the B-spline-like basis functions, the polynomial framework of $\C{1}$ cubic Powell--Sabin splines can be readily generalized to a rational one \cite{rps3_groselj_25}. For each $v_i \in \vertices$ and $r \in \set{1,2,3}$, let us associate with $B_{i,r}^v$ a weight $w_{i,r}^v \in (0, \infty)$, and, for each $e_l \in \edges$ and $r \in \set{1,2}$, let us associate with $B_{l,r}^e$ a weight $w_{l,r}^e \in (0, \infty)$. Moreover, let the spline $W \in \mathbb{S}_3(\PST)$ be given by
\begin{equation*}
W = \sum_{i=1}^{|\vertices|} \sum_{r=1}^3 w_{i,r}^v B_{i,r}^v + \sum_{l=1}^{|\edges|} \sum_{r=1}^2 w_{l,r}^e B_{l,r}^e.
\end{equation*}
The spline $W$ is positive everywhere. Based on this spline and the weights, we define a set of rational basis functions. Namely,
\begin{itemize}
\item three basis functions $N_{i,r}^v: \Theta \rightarrow \R$, $r = 1, 2, 3$, for each $v_i \in \vertices$, given by
\begin{equation*}
N_{i,r}^v = \frac{w_{i,r}^v}{W} B_{i,r}^v,
\end{equation*}

\item two basis functions $N_{l,r}^e: \Theta \rightarrow \R$, $r = 1, 2$, for each $e_l \in \edges$, given by
\begin{equation*}
N_{l,r}^e = \frac{w_{l,r}^e}{W} B_{l,r}^e.
\end{equation*}
\end{itemize}
By the assumptions on the weights, these functions are well-defined and nonnegative. The functions are rational and piecewise cubic in numerator and denominator. If all $w_{i,r}^v$ and all $w_{l,r}^e$ are equal, $N_{i,r}^v = B_{i,r}^v$ and $N_{l,r}^e = B_{l,r}^e$. For a fixed set of weights $w_{i,r}^v$ and $w_{l,r}^e$, we denote by $\mathbb{S}_3^w(\PST)$ the vector space spanned by $N_{i,r}^v$ and $N_{l,r}^e$.

Finally, let us note that the gradients of $N_{i,r}^v$ and $N_{l,r}^e$ can be computed based on the gradients of $B_{i,r}^v$ and $B_{l,r}^e$ in a straightforward manner. Namely,
\begin{equation*}
\grad{N_{i,r}^v} = \frac{w_{i,r}^v}{W} \grad{B_{i,r}^v} - \frac{N_{i,r}^v}{W} \, \grad{W}, 
\qquad
\grad{N_{l,r}^e} = \frac{w_{l,r}^e}{W} \grad{B_{l,r}^e} - \frac{N_{l,r}^e}{W} \, \grad{W}.
\end{equation*}

\section{Prelude to approximation}
\label{sec:approx}

\subsection{Geometry mapping}
\label{sec:geometry}

Let $\mathbb{S}_3^w(\PST)^3 = \mathbb{S}_3^w(\PST) \times \mathbb{S}_3^w(\PST) \times \mathbb{S}_3^w(\PST)$. The basis functions $N_{i,r}^v$ and $N_{l,r}^e$ provide a flexible tool to define a geometry mapping $F \in \mathbb{S}_3^w(\PST)^3$ such that $\Omega = F(\Theta) \subset \R^3$ is a bounded surface with connected boundary $\partial \Omega = F(\partial \Theta)$. Suppose $F$ is a diffeomorphism given by
\begin{equation}
\label{eq:geometry_mapping}
F = \sum_{i=1}^{|\vertices|} \sum_{r=1}^3 (x_{i,r}^v, y_{i,r}^v, z_{i,r}^v) N_{i,r}^v + \sum_{l=1}^{|\edges|} \sum_{r=1}^2 (x_{l,r}^e, y_{l,r}^e, z_{l,r}^e) N_{l,r}^e
\end{equation}
for some $x_{i,r}^v, y_{i,r}^v, z_{i,r}^v \in \R$ and $x_{l,r}^e, y_{l,r}^e, z_{l,r}^e \in \R$. Let $J_F: \Theta \rightarrow \R^{3 \times 2}$ denote the mapping assigning to each point the Jacobian matrix of $F$, which is assumed to be of full rank. Then, the mapping $K_F: \Theta \rightarrow \R^{2 \times 2}$ given by $K_F = {J_F}^T J_F$ assigns to each point a symmetric positive definite matrix, and the function $\kappa_F: \Theta \rightarrow \R$ given by $\kappa_F = \sqrt{\det(K_F)}$ is well-defined and positive.

The surface representation in \eqref{eq:geometry_mapping} is well-suited for freeform modeling. As the basis functions form a convex partition of unity on $\Theta$, one can consider $(x_{i,r}^v, y_{i,r}^v, z_{i,r}^v)$ and $(x_{l,r}^e, y_{l,r}^e, z_{l,r}^e)$ in \eqref{eq:geometry_mapping} as control points to modify the surface. Additionally, the weights appearing in $N_{i,r}^v$ and $N_{l,r}^e$ can be used to finetune the shape.

In the special case where all $z_{i,r}^v$ and all $z_{l,r}^e$ are zero, we can consider $F$ as a mapping to $\R^2$, \ie{}, $F \in \mathbb{S}_3^w(\PST)^2$, and $\Omega$ as a bounded domain in $\R^2$. Under the previous assumptions, the mapping $J_F$ then assigns to each point from $\Theta$ an invertible matrix, and the function $\det(J_F): \Theta \rightarrow \R$ is either positive or negative. Commonly, the former is required for $F$ to preserve the orientation.

Often, one describes a planar domain $\Omega$ by specifying the parametrization of its boundary in the form of a bijective spline curve $C: \partial \Theta \rightarrow \partial \Omega$ and defines $F$ on $\Theta$ by minimizing a suitable energy functional under the conditions that $F(\partial \Theta) = C(\partial \Theta) = \partial \Omega$ and $\det(J_F)$ is positive. This is a sufficient condition for the bijectivity of $F$ \cite{kestelman_71}. The following example is taken from \cite{ps2_speleers_15}, where quadratic Powell--Sabin splines are used, and adapted to the cubic case.

\begin{figure}[t]
\centering

\begin{subfigure}[b]{0.32\textwidth}
\centering
\includegraphics[width=\textwidth]{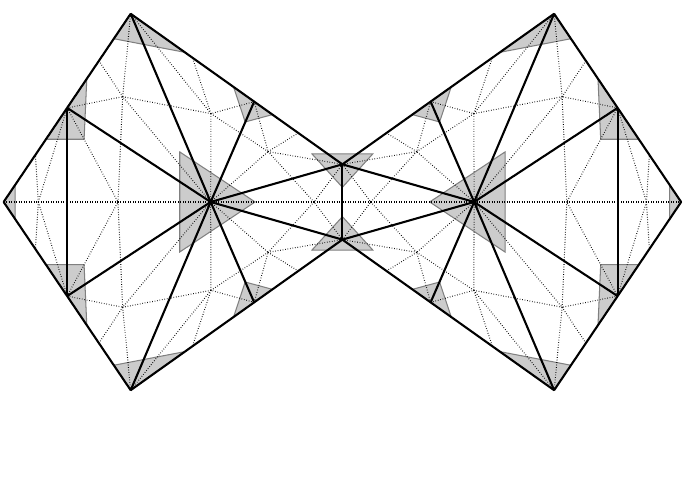}
\caption{The domain $\Theta$ and its configuration.}
\label{fig:puzzle_tri}
\end{subfigure}
\hfill
\begin{subfigure}[b]{0.32\textwidth}
\centering
\includegraphics[width=\textwidth]{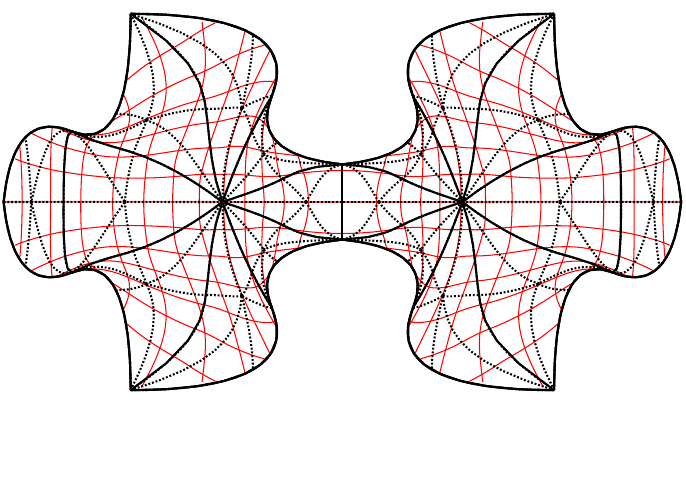}
\caption{The domain $\Omega$ parametrized by $F$.}
\label{fig:puzzle_curvtri}
\end{subfigure}
\hfill
\begin{subfigure}[b]{0.32\textwidth}
\centering
\includegraphics[width=\textwidth]{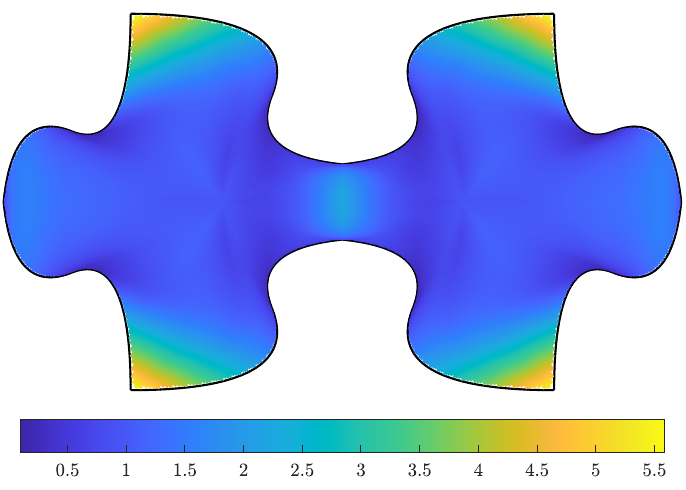}
\caption{The values of $\det(J_F)$ depicted on $\Omega$.}
\label{fig:puzzle_jacobi}
\end{subfigure}

\caption{A cubic Powell--Sabin spline parametrization $F: \Theta \rightarrow \Omega$ described in Example~\ref{ex:puzzle}.}
\label{fig:puzzle}
\end{figure}

\begin{example}
\label{ex:puzzle}
The domain $\Omega \subset \R^2$ shown in Figure~\ref{fig:puzzle} is specified by quadratic B\'{e}zier curves given in \cite[Table 1]{ps2_speleers_15}. The mapping $F = (F_1, F_2)$ of the form \eqref{eq:geometry_mapping} is defined on the parametric domain $\Theta \subset \R^2$ shown in Figure~\ref{fig:puzzle_tri} together with the triangulation $\triangle$ (solid lines), the Powell--Sabin triangulation $\PST$ (dotted lines), and a choice of triangles used in the construction of the Powell--Sabin B-splines (gray triangles associated with the vertices of $\triangle$). The weights $w_{i,r}^v$ and $w_{l,r}^e$ of $N_{i,r}^v$ and $N_{l,r}^e$ are all chosen to be equal to one, and a subset of control points of $F$ is fixed to exactly reproduce the boundary. The remaining control points are obtained by minimizing the normalized Winslow functional
\begin{equation*}
\omega(F) = \frac{1}{\mathrm{area}(\Theta)} \int_\Theta \frac{\tnorm{\nabla F_1}^2 + \tnorm{\nabla F_2}^2}{\abs{\det(J_F)}} \, \mathrm{d}\Theta,
\end{equation*}
subjected to the condition that $\det(J_F)$ is positive. The resulting parametrization is depicted in Figure~\ref{fig:puzzle_curvtri} by isoparametric lines in red color and mapped edges of the Powell--Sabin triangulation in black color. The attained value of $\omega(F)$ is $2.5827$ (in general, the value of $\omega(F)$ cannot be smaller than $2$) and the values of $\det(J_F)$ range in the interval $[0.1, 5.6]$ as shown in Figure~\ref{fig:puzzle_jacobi}. These results indicate that $F$ is a parametrization of good quality.
\end{example}

If $G: \Theta \rightarrow \Omega$ is a known regular parametrization of the domain, one can use a projector to $\mathbb{S}_3(\PST)$ to define $F \in \mathbb{S}_3^w(\PST)^3$. Suitable techniques for the construction of projection operators were investigated in \cite{ps3_groselj_18} but have to be applied with caution to produce $F$ that meets the assumptions and is of good quality. In the special case where $G \in \mathbb{S}_3^w(\PST)^3$, the use of a projector gives rise to a convenient procedure to transform $G = F$ into the form \eqref{eq:geometry_mapping}; see Examples \ref{ex:annulus} and \ref{ex:cylinder}. Otherwise the image of $F$ is only an approximation of $\Omega$. However, for planar domains, if the boundary of $\Omega$ is a $\C{1}$ cubic spline curve, it is possible to achieve that $\Omega = F(\Theta) = G(\Theta)$ although $F \neq G$; see Example~\ref{ex:pentagon}.

\begin{figure}[t]
\centering

\begin{subfigure}[b]{0.32\textwidth}
\centering
\includegraphics[width=\textwidth]{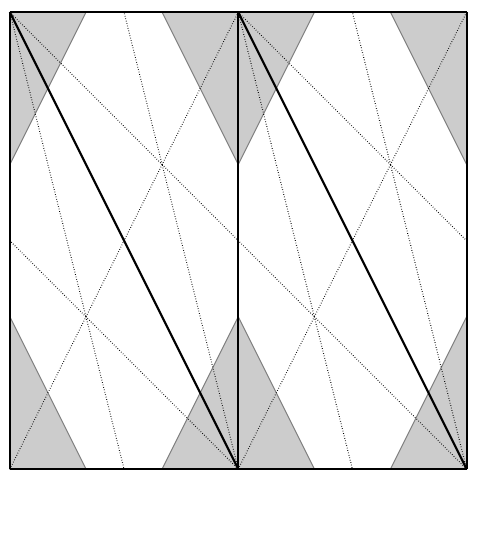}
\caption{The domain $\Theta$ and its configuration.}
\label{fig:annulus_tri}
\end{subfigure}
\hfill
\begin{subfigure}[b]{0.32\textwidth}
\centering
\includegraphics[width=\textwidth]{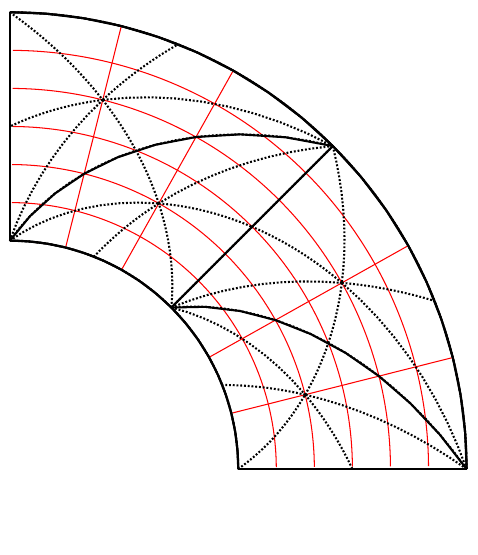}
\caption{The domain $\Omega$ parametrized by $F$.}
\label{fig:annulus_curvtri}
\end{subfigure}
\hfill
\begin{subfigure}[b]{0.32\textwidth}
\centering
\includegraphics[width=\textwidth]{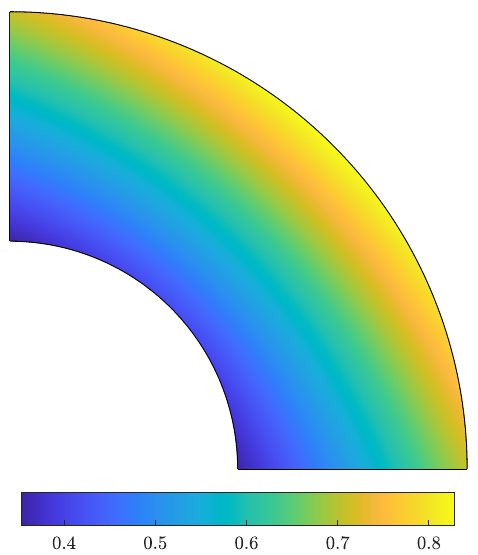}
\caption{The values of $\kappa_F$ depicted on $\Omega$.}
\label{fig:annulus_jacobi}
\end{subfigure}

\caption{A cubic Powell--Sabin spline parametrization $F: \Theta \rightarrow \Omega$ described in Example~\ref{ex:annulus}.}
\label{fig:annulus}
\end{figure}

\begin{example}
\label{ex:annulus}
Consider the annulus bounded by two circles with centers at $(0,0)$ and radii $r_1$ and $r_2$, $0 < r_1 < r_2$. Let $\Omega \subset \R^2$ be the quarter of this annulus positioned in the first quadrant of the coordinate system, \ie{},
\begin{equation*}
\Omega = \defset{(x,y) \in [0, r_2] \times [0, r_2]}{r_1 \leq \sqrt{x^2+y^2} \leq r_2}.
\end{equation*}
To parametrize $\Omega$, we consider a specific quadratic NURBS curve $C: [0, 1] \rightarrow \R^2$ defined on the knot set $\set{0, 0, 0, \frac{1}{2}, 1, 1, 1}$ and specified by the control points $(1,0)$, $(1, \sqrt{2}-1)$, $(\sqrt{2} - 1, 1)$, $(0, 1)$ and the weights $1$, $\frac{1}{4} (2 + \sqrt{2})$, $\frac{1}{4} (2 + \sqrt{2})$, $1$. It is straightforward to verify that this is an exact parametrization of the unit circular arc in the first quadrant with the center at $(0,0)$. Suppose $C_1 = r_1 C$ and $C_2 = r_2 C$. Clearly, $\Omega$ can be parametrized as a tensor-product NURBS surface $F: [0, 1] \times [0, 1] \rightarrow \R^2$ of degree $2 \times 1$ ruled by $C_1$ and $C_2$. Since $F$ is $\C{1}$ continuous and piecewise cubic polynomial in numerator and denominator, it can be represented in the form \eqref{eq:geometry_mapping} by the procedure provided in \cite{rps3_groselj_25}. More precisely, the parametric domain $\Theta = [0,1] \times [0,1]$ is triangulated by four triangles as shown in Figure~\ref{fig:annulus_tri}, and the parametrization of $\Omega$ for $r_1 = \frac{1}{2}$ and $r_2 = 1$ is shown in Figure~\ref{fig:annulus_curvtri}.
\end{example}

\begin{figure}[t]
\centering

\begin{subfigure}[b]{0.32\textwidth}
\centering
\includegraphics[width=\textwidth]{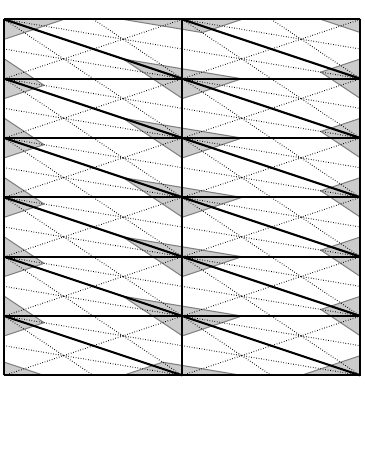}
\caption{The domain $\Theta$ and its configuration.}
\label{fig:cylinder_tri}
\end{subfigure}
\hfill
\begin{subfigure}[b]{0.32\textwidth}
\centering
\includegraphics[width=\textwidth]{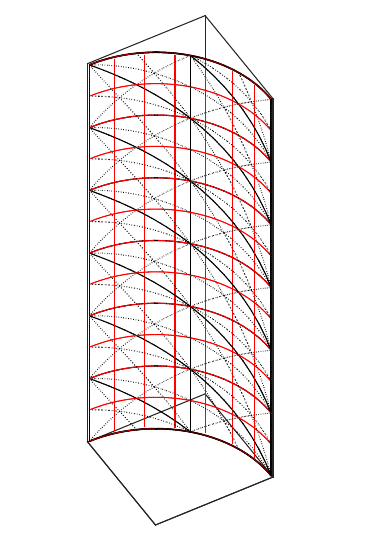}
\caption{The domain $\Omega$ parametrized by $F$.}
\label{fig:cylinder_curvtri}
\end{subfigure}
\hfill
\begin{subfigure}[b]{0.32\textwidth}
\centering
\includegraphics[width=\textwidth]{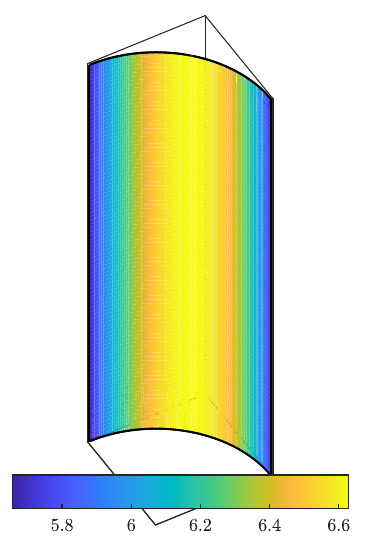}
\caption{The values of $\kappa_F$ depicted on $\Omega$.}
\label{fig:cylinder_jacobi}
\end{subfigure}

\caption{A cubic Powell--Sabin spline parametrization $F: \Theta \rightarrow \Omega$ described in Example~\ref{ex:cylinder}.}
\label{fig:cylinder}
\end{figure}

\begin{example}
\label{ex:cylinder}
Consider the cylinder with radius $r > 0$ and height $h > 0$ that is revolved around the $z$-axis. Let $\Omega \subset \R^3$ be the quarter of this cylinder positioned in the first octant of the coordinate system, \ie{},
\begin{equation*}
\Omega = \defset{(x,y,z) \in [0, r] \times [0, r] \times [0, h]}{\sqrt{x^2 + y^2} = r}.
\end{equation*}
We can parametrize $\Omega$ as a ruled surface $F: [0, 1] \times [0,1] \rightarrow \R^3$ bounded by two spatial NURBS curves $C_1, C_2: [0, 1] \rightarrow \R^3$. With the help of the curve $C$ from Example~\ref{ex:annulus}, these two curves can be specified as $C_1 = (r C, h)$ and $C_2 = (r C, 0)$. Then, using the same procedure as described in Example~\ref{ex:annulus}, the parametrization $F$ of $\Omega$ can be exactly represented in the form \eqref{eq:geometry_mapping}. For a particular case, the parametric domain $\Theta = [0,1] \times [0,1]$ is triangulated by $12$ triangles as shown in Figure~\ref{fig:cylinder_tri}, and the parametrization of $\Omega$ for $r = 1$ and $h = 4$ is shown in Figure~\ref{fig:cylinder_curvtri}.
\end{example}

\begin{figure}[t]
\centering

\begin{subfigure}[b]{0.31\textwidth}
\centering
\includegraphics[width=\textwidth]{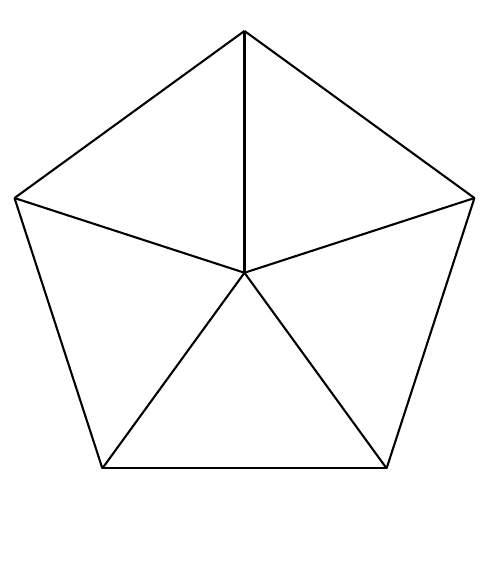}
\caption{The domain $\Theta$ and its triangulation.}
\label{fig:pentagon_bb_tri}
\end{subfigure}
\hfill
\begin{subfigure}[b]{0.31\textwidth}
\centering
\includegraphics[width=\textwidth]{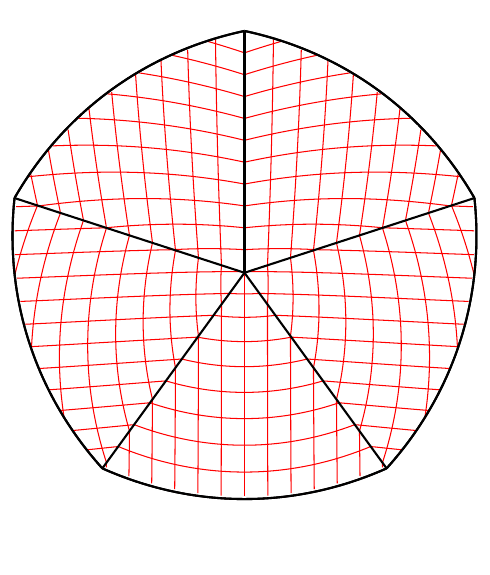}
\caption{The domain $\Omega$ as parametrized by $G$.}
\label{fig:pentagon_bb_curvtri}
\end{subfigure}
\hfill
\begin{subfigure}[b]{0.31\textwidth}
\centering
\includegraphics[width=\textwidth]{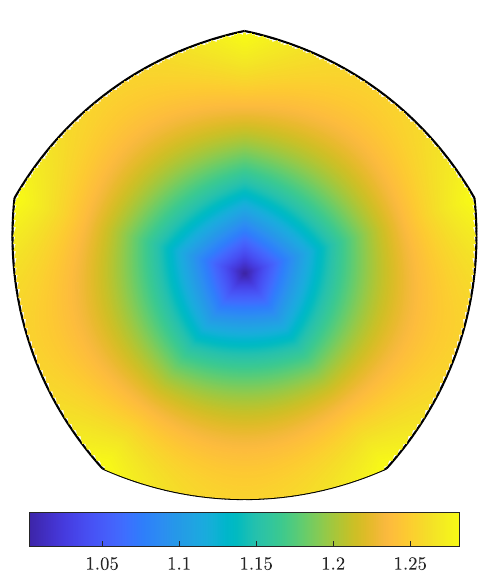}
\caption{The values of $\det(J_G)$ depicted on $\Omega$.}
\label{fig:pentagon_bb_jacobi}
\end{subfigure}

\begin{subfigure}[b]{0.31\textwidth}
\centering
\includegraphics[width=\textwidth]{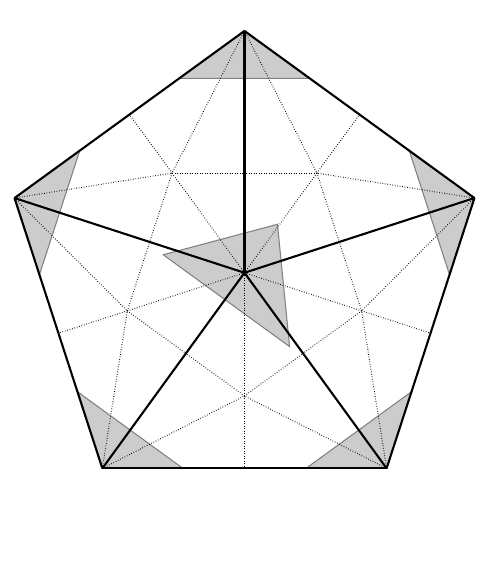}
\caption{The domain $\Theta$ and its configuration.}
\label{fig:pentagon_ps_tri}
\end{subfigure}
\hfill
\begin{subfigure}[b]{0.31\textwidth}
\centering
\includegraphics[width=\textwidth]{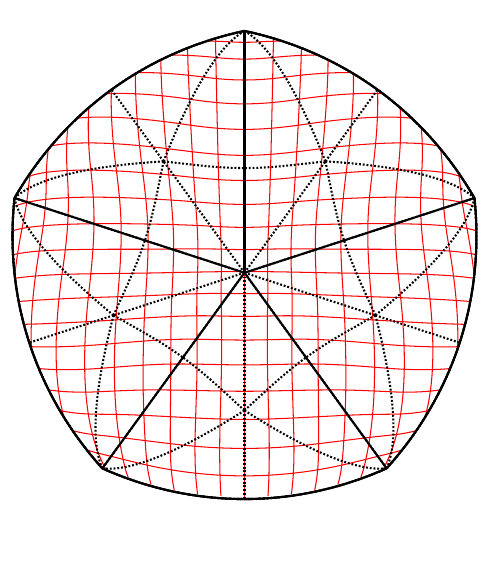}
\caption{The domain $\Omega$ as parametrized by $F$.}
\label{fig:pentagon_ps_curvtri}
\end{subfigure}
\hfill
\begin{subfigure}[b]{0.31\textwidth}
\centering
\includegraphics[width=\textwidth]{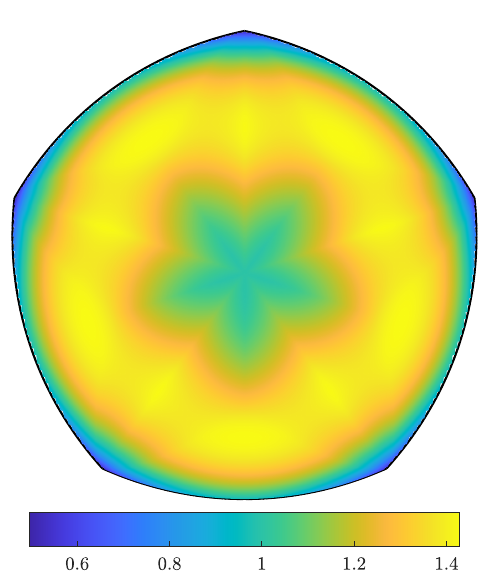}
\caption{The values of $\det(J_F)$ depicted on $\Omega$.}
\label{fig:pentagon_ps_jacobi}
\end{subfigure}

\caption{The $\C{0}$ cubic spline parametrization $G: \Theta \rightarrow \Omega$ (top) and the $\C{1}$ cubic Powell--Sabin spline parametrization $F: \Theta \rightarrow \Omega$ (bottom) that are described in Example~\ref{ex:pentagon}.}
\label{fig:pentagon}
\end{figure}

\begin{example}
\label{ex:pentagon}
Let us consider the cubic triangular B\'{e}zier patch in $\R^2$ defined on the domain triangle with the vertices $(0,0)$, $(-\sin(\frac{1}{5} \pi), -\cos(\frac{1}{5} \pi))$, $(\sin(\frac{1}{5} \pi), -\cos(\frac{1}{5} \pi))$ and specified by the control points
\begin{equation*}
\tfrac{j}{3}\bracket{-\sin(\tfrac{j - 2k}{5j} \pi), -\cos(\tfrac{j - 2k}{5j} \pi)}, \qquad j = 0, 1, 2, 3, \quad k = 0, 1, \ldots, j.
\end{equation*}
Let $\Omega \subset \R^2$ be the union of this B\'{e}zier patch and its four consecutive rotations by the angle $\frac{2}{5} \pi$ around $(0,0)$. Similarly, let $\Theta \subset \R^2$ be the union of the domain triangle and its four consecutive rotations. The parametric domain $\Theta$ with the induced triangulation and the domain $\Omega$ are shown in Figure~\ref{fig:pentagon} (top). By this construction, $\Omega$ is parametrized by a $\C{0}$ piecewise cubic mapping $G: \Theta \rightarrow \Omega$ with $\det(J_G) \in [1.0, 1.3]$ and $\omega(G) = 2.0666$. In order to obtain a $\C{1}$ parametrization $F: \Theta \rightarrow \Omega$ of the form \eqref{eq:geometry_mapping} with all weights equal to one, we follow \cite[Recipe 2]{ps3_groselj_18} to set up a quasi-interpolation operator that projects $G$ to $\mathbb{S}_3(\PST)^2$ and satisfies $F(\partial \Theta) = G(\partial \Theta)$. The resulting mapping $F$ is of the form \eqref{eq:geometry_mapping} with $\det(J_F) \in [0.5, 1.5]$ and $\omega(F) = 2.0668$. Hence $F$ is a $\C{1}$ bijection and $F(\Theta) = G(\Theta) = \Omega$. The parametrization $F$ is depicted in Figure~\ref{fig:pentagon} (bottom).
\end{example}

\subsection{Approximation function and boundary conditions}
\label{sec:approx_boundary}

To define an approximation function on a domain $\Omega \subset \R^3$ parametrized by a geometry mapping $F \in \mathbb{S}_3^w(\PST)^3$, we follow the isoparametric approach and set a function $\varphi: \Omega \rightarrow \R$ as $\varphi = s \circ F^{-1}$ for some $s \in \mathbb{S}_3^w(\PST)$. We consider $s$ as
\begin{equation}
\label{eq:spline_approx}
s = \sum_{i = 1}^{|\vertices|} \sum_{r = 1}^3 c_{i,r}^v N_{i,r}^v + \sum_{l = 1}^{|\edges|} \sum_{r = 1}^2 c_{l,r}^e N_{l,r}^e,
\end{equation}
where the coefficients $c_{i,r}^v \in \R$ and $c_{l,r}^e \in \R$ are free parameters in the approximation process. Typically, a subset of these coefficients is fixed by boundary conditions. More particularly, in Section~\ref{sec:problems} we deal with two types of boundary conditions. The first is the standard Dirichlet condition $\varphi = g_0$ on $\partial \Omega$ for some function $g_0: \partial \Omega \rightarrow \R$, which implies
\begin{subequations}
\label{eq:boundary_conditions}
\begin{equation}
\label{eq:boundary_dirichlet}
s = g_0 \circ F \qquad \text{on $\partial \Theta.$}
\end{equation}
The second is the Dirichlet condition supplemented by $\nabla_\Omega \varphi \cdot n_{\partial \Omega} = g_1$ on $\partial \Omega$ for some function $g_1: \partial \Omega \rightarrow \R$, where $\nabla_\Omega$ is the surface gradient \cite{Delfour2011} and $n_{\partial \Omega}: \partial \Omega \rightarrow \R^3$ is the vector field specifying the outward unit vector normal to the boundary $\partial \Omega$, \ie{}, the vector collinear to the cross product of the tangent to the boundary and the normal to the surface $\Omega$ at the considered point on $\partial \Omega$. In relation to the parametrization $F$ of $\Omega$, the surface gradient can be expressed as
\begin{equation*}
\nabla_\Omega \varphi \circ F = \nabla_\Omega (s \circ F^{-1}) \circ F = J_F {K_F}^{-1} \nabla s.
\end{equation*}
Moreover, let us denote by $\tau: \partial \Theta \rightarrow \R^2$ the vector field specifying the tangent and by $\nu: \partial \Theta \rightarrow \R^2$ the vector field specifying the outward unit normal to the piecewise linear boundary of $\Theta$. Since $J_F \tau$ is tangent to the boundary of $\Omega$, there exists a scalar field $\gamma: \partial \Theta \rightarrow \R$ such that $n_{\partial \Omega} \circ F = \gamma J_F {K_F}^{-1} \nu$. This implies
\begin{equation}
\label{eq:boundary_normal}
\nabla s \cdot \gamma {K_F}^{-1} \nu = g_1 \circ F \qquad \text{on $\partial \Theta.$}
\end{equation}
\end{subequations}
Thus, in order to identify the coefficients affected by the conditions in \eqref{eq:boundary_conditions}, it is important to have a better understanding of the supports of the basis functions $N_{i,r}^v$ and $N_{l,r}^e$ and their gradients.

It is clear from the definition of $N_{i,r}^v$ and $N_{l,r}^e$ that their supports match the supports of the corresponding basis functions $B_{i,r}^v$ and $B_{l,r}^e$. Thus, following the construction in \cite[Section 3]{ps3_speleers_15}, the support of $N_{i,r}^v$ is contained in the union of triangles of $\triangle$ that have the vertex at $v_i$ (see Figure~\ref{fig:psbv}), and the support of $N_{l,r}^e$ is contained in the union of triangles of $\triangle$ that have the edge $e_l$ in common (see Figure~\ref{fig:psbe}). Moreover, at a vertex $v_i \in \vertices$ only the basis functions $N_{i,r}^v$ have nonzero values and gradients, and on an edge $e_l = \conv{v_i, v_j} \in \edges$ only the basis functions $N_{i,r}^v$, $N_{j,r}^v$, and  $N_{l,r}^e$ may have nonzero values and gradients. If $e_l$ is a boundary edge, then one of $N_{l,r}^e$, $r = 1, 2$, is zero on $e_l$ (see Figure~\ref{fig:psbe2}), while the values of $N_{i,r}^v$ and $N_{j,r}^v$ on $e_l$ depend on the choice of triangles that define the configuration of $B_{i,r}^v$ and $B_{l,r}^e$.

The configuration of the basis functions is determined by a set of triangles in $\R^2$, one for each vertex $v_i \in \vertices$ (see Figure~\ref{fig:pstri}). The triangle $q_i = \conv{q_{i,1}, q_{i,2}, q_{i,3}} \subset \R^2$ associated with $v_i$ contains the vertex $v_i$ and, for each $t_m = \conv{v_i, v_l^e, v_m^t} \in \PST$ with a vertex at $v_i$, also the points $\frac{2}{3} v_i + \frac{1}{3} v_l^e$ and $\frac{2}{3} v_i + \frac{1}{3} v_m^t$. This ensures the nonnegativity of the basis functions \cite[Theorem 3]{ps3_speleers_15} but also gives some further insight into the values of $B_{i,r}^v$ on the edge $e_l$. From \cite[Sections 3.2 and 3.3]{ps3_speleers_15} it follows that they are completely specified by two barycentric coordinates with respect to the triangle $q_i$, \ie{}, the barycentric coordinate of the point $v_i$ and the barycentric coordinate of the point $\frac{2}{3} v_i + \frac{1}{3} v_l^e$, both corresponding to the vertex $q_{i,r}$. In the specific case where the triangle $q_i$ is positioned in such a way that its edge opposite $q_{i,r}$ lies on the line containing the edge $e_l$, both these coordinates are equal to zero and hence the restriction of $B_{i,r}^v$ to $e_l$ is zero. Due to the constraints on the triangle $q_i$ this is possible only when $v_i$ is a vertex on the boundary of $\Theta$. Moreover, if $B_{i,r}^v$ is to be equal to zero on the entire boundary, the two boundary edges of $\triangle$ that meet at $v_i$ must lie on the same line. In this situation only two of the three basis functions associated with $v_i$ are nonzero on $\partial \Theta$. In order to minimize the number of coefficients affected by the boundary conditions, we adhere to such positioning of triangles on the boundary. For an illustration, see Figure~\ref{fig:pstri} and compare Figures \ref{fig:psbv2} and \ref{fig:psbv3}.

To summarize, let the set $\vertices^b$ of boundary vertices of $\triangle$ be considered as the disjoint union of $\widehat{\vertices}^b$ and $\bar{\vertices}^b$, where $v_i \in \vertices^b$ belongs to $\widehat{\vertices}^b$ if the two boundary edges of $\triangle$ that meet at $v_i$ do not lie on the same line, and to $\bar{\vertices}^b$ otherwise. Then, a suitable configuration of triangles $q_i$ ensures that the number of nonzero basis functions on $\partial \Theta$ is $b_0 = 3 |\widehat{\vertices}^b| + 2 |\bar{\vertices}^b| + |\bedges| = 4|\widehat{\vertices}^b| + 3|\bar{\vertices}^b|$. The number of basis functions that are nonzero on $\partial \Theta$ or have a nonvanishing gradient on $\partial \Theta$ is $b_1 = 3|\vertices^b| + 2|\bedges| = 5|\vertices^b|$, regardless of the configuration.

In general, the boundary conditions in \eqref{eq:boundary_conditions} cannot be satisfied exactly. If the condition \eqref{eq:boundary_dirichlet} is imposed, we utilize the least squares fitting on a dense set of points on $\partial \Theta$ to determine the identified $b_0$ coefficients and compute $s_0 \in \mathbb{S}_3^w(\PST)$ that approximates $g_0 \circ F$. Similarly, if both \eqref{eq:boundary_dirichlet} and \eqref{eq:boundary_normal} are imposed, we determine the identified $b_1$ coefficients and compute $s_1 \in \mathbb{S}_3^w(\PST)$ such that $s_1$ and $\nabla s_1$ provide the best least squares fit to $g_0 \circ F$ and $g_1 \circ F$ in the sense of \eqref{eq:boundary_conditions}.

To simplify the notation in the rest of the paper, we rename and reindex the basis functions $N_{i,r}^v$ and $N_{l,r}^e$ to $N_k$, $k = 1, 2, \ldots, 3|\vertices| + 2 |\edges|$, in such a way that when \eqref{eq:boundary_dirichlet} is imposed, the first $n_0 = 3|\vertices| + 2 |\edges| - b_0$ basis functions are zero on $\partial \Theta$, and when \eqref{eq:boundary_dirichlet} and \eqref{eq:boundary_normal} are imposed, the first $n_1 = 3|\vertices| + 2 |\edges| - b_1$ basis functions are zero on $\partial \Theta$ and their gradients are zero on $\partial \Theta$ as well. Thus \eqref{eq:spline_approx} can be rewritten as
\begin{equation}
\label{eq:spline_approx_boundary}
s = s_d + \sum_{k = 1}^{n_d} c_k N_k, \qquad d = 0, 1,
\end{equation}
and the conditions $c_k \in \R$ represent the remaining degrees of freedom.

\section{Approximation problems}
\label{sec:problems}

\subsection{Poisson problems}

Let $\Delta_\Omega = \nabla_\Omega \cdot \nabla_\Omega$ be the Laplace--Beltrami operator on the surface domain $\Omega$ parametrized by the geometry mapping $F \in \mathbb{S}_3^w(\PST)^3$. We are interested in finding an approximation $s \in \mathbb{S}_3^w(\PST)$ to a function $u: \Omega \rightarrow \R$ satisfying
\begin{equation}
\label{eq:poisson_problem}
\left\{
\begin{aligned}
-\Delta_\Omega u &= f && \text{on} \ \Omega \setminus \partial \Omega, \\
u &= g_0 && \text{on} \ \partial \Omega,
\end{aligned}
\right.
\end{equation}
for given functions $f: \Omega \rightarrow \R$ and $g_0: \partial \Omega \rightarrow \R$. To this end, we follow \cite{iga_dede_15} and consider the weak formulation of \eqref{eq:poisson_problem} under suitable smoothness conditions on $f$ and $g_0$, in particular $f \in L^2(\Omega)$. After the discretization based on $\mathbb{S}_3^w(\PST)$ and the pull-back to the parametric domain $\Theta$, we arrive at the conditions
\begin{equation}
\label{eq:isofem_conditions}
\int_\Theta \grad{s} \cdot \left({K_F}^{-1} \, \grad{N_j} \right) \, \kappa_F \, \mathrm{d}\Theta = \int_\Theta (f \circ F) \, N_j \, \kappa_F \, \mathrm{d}\Theta, \qquad j = 1, 2, \ldots, n_0,
\end{equation}
and the boundary conditions \eqref{eq:boundary_dirichlet} imposed on the spline $s$ expressed in the form \eqref{eq:spline_approx_boundary} for $d = 0$. Then, by \eqref{eq:isofem_conditions}, the coefficients of $s$ collected in the vector $c = (c_k)_k \in \R^{n_0}$ are the solution of the system of linear equations $A c = b$, where the entries of the matrix $A = (a_{j,k})_{j,k} \in \R^{n_0 \times n_0}$ are given by
\begin{equation*}
a_{j,k} = \int_\Theta \grad{N_k} \cdot \left({K_F}^{-1} \, \grad{N_j}\right) \, \kappa_F \, \mathrm{d}\Theta,
\end{equation*}
and the entries of the vector $b = (b_j)_j \in \R^{n_0}$ are given by
\begin{equation*}
b_j = \int_\Theta \bracket{(f \circ F) \, N_j - \grad{s_0} \cdot \left({K_F}^{-1} \, \grad{N_j}\right)} \, \kappa_F \, \mathrm{d}\Theta.
\end{equation*}

When $\Omega = \Theta$ and $F: \Theta \rightarrow \Omega$ is the identity mapping in $\R^2$, the weak form discretization \eqref{eq:isofem_conditions} reduces to the standard Ritz--Galerkin method on $\mathbb{S}_3(\PST)$. In the following example we compare this method with the common approach based on cubic Lagrange elements. We note that the cubic Lagrange elements over $\triangle$ have a very comparable number of degrees of freedom, \ie{}, $n_0 \approx (|\vertices| - |\bvertices|) + 2(|\edges| - |\bedges|) + |\triangles|$ because $|\triangles| = 2|\vertices| - |\bvertices| - 2$. The results indicate that the Powell--Sabin elements provide a slight but consistent improvement over the Lagrange elements.

\begin{figure}[t]
\centering

\begin{subfigure}[b]{0.19\textwidth}
\centering
\includegraphics[width=\textwidth]{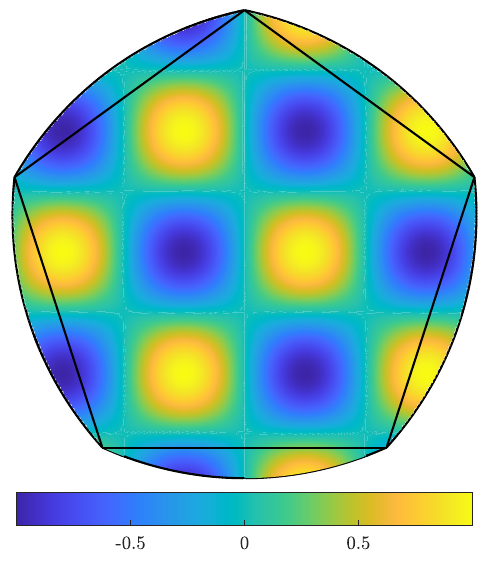}
\caption{Examples \ref{ex:fem} and \ref{ex:isofem}.}
\label{fig:pentagon_fun}
\end{subfigure}
\hfill
\begin{subfigure}[b]{0.195\textwidth} 
\centering
\includegraphics[width=\textwidth]{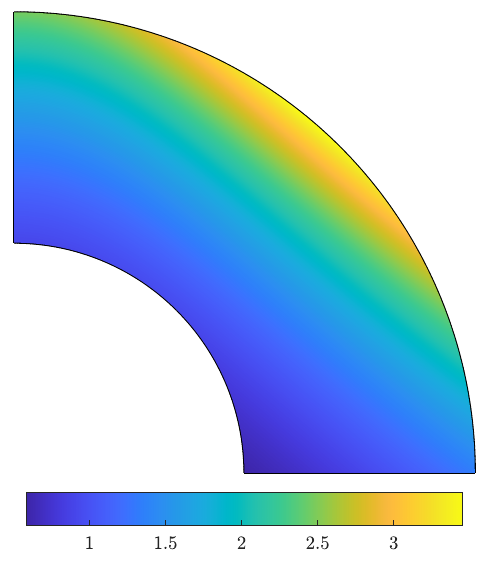}
\caption{Examples \ref{ex:poi2} and \ref{ex:bih2}.}
\label{fig:poi2_fun}
\end{subfigure}
\hfill
\begin{subfigure}[b]{0.19\textwidth}
\centering
\includegraphics[width=\textwidth]{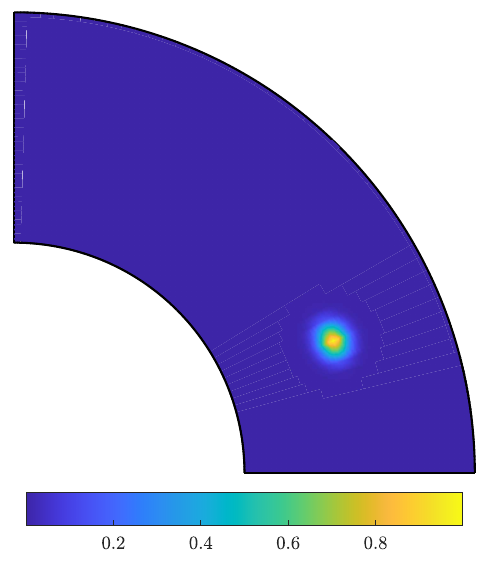}
\caption{Example~\ref{ex:poi2_locref}.}
\label{fig:poi2_locref_fun}
\end{subfigure}
\hfill
\begin{subfigure}[b]{0.19\textwidth}
\centering
\includegraphics[width=\textwidth]{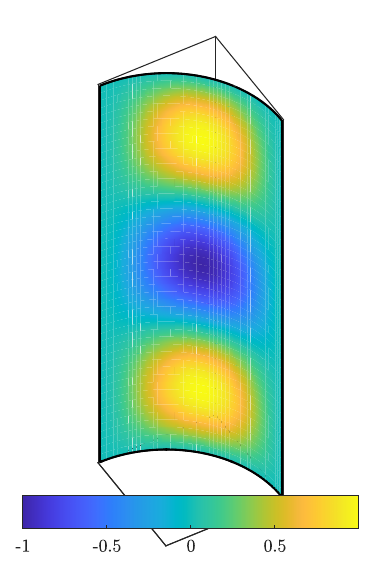}
\caption{Example~\ref{ex:poi3}.}
\label{fig:poi3_fun}
\end{subfigure}
\hfill
\begin{subfigure}[b]{0.19\textwidth}
\centering
\includegraphics[width=\textwidth]{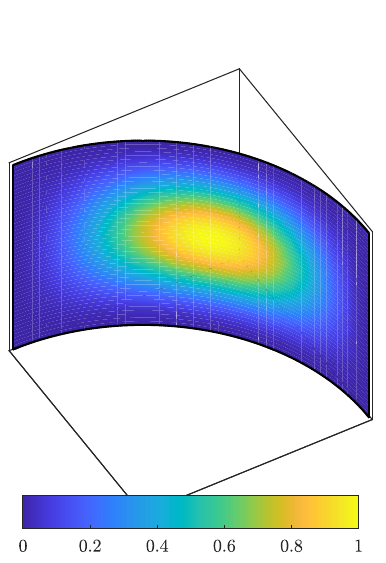}
\caption{Example~\ref{ex:bih3}.}
\label{fig:bih3_fun}
\end{subfigure}

\caption{Graphs of the functions $u$ used in Examples \ref{ex:fem}--\ref{ex:bih3}.}
\label{fig:fun}
\end{figure}

\begin{figure}[ht]
\centering

\begin{subfigure}[b]{0.48\textwidth}
\centering
\includegraphics[width=\textwidth]{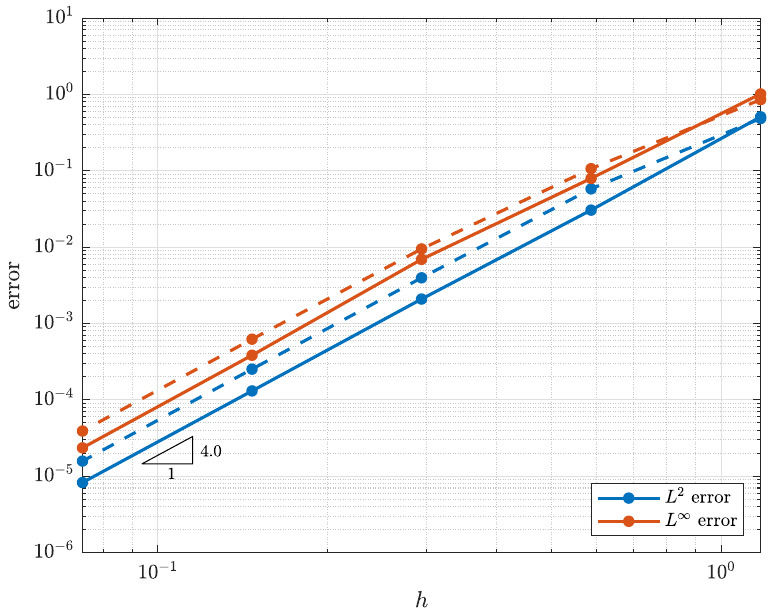}
\caption{Example \ref{ex:fem}.}
\label{fig:fem_err}
\end{subfigure}
\hfill
\begin{subfigure}[b]{0.48\textwidth}
\centering
\includegraphics[width=\textwidth]{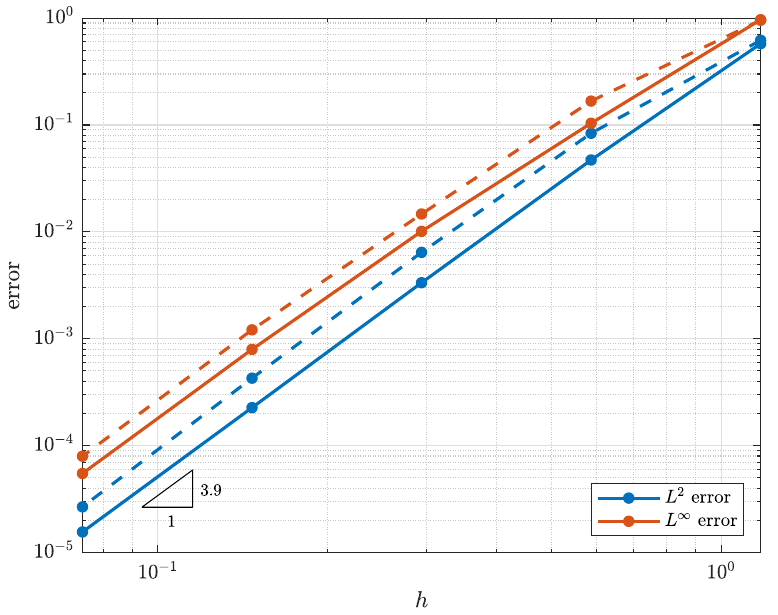}
\caption{Example \ref{ex:isofem}.}
\label{fig:isofem_err}
\end{subfigure}

\caption{Approximation errors for the problems described in Example~\ref{ex:fem} (left) and Example~\ref{ex:isofem} (right). The errors are presented with respect to the mesh size $h$. The solid and dashed lines correspond to the results obtained by the Powell--Sabin and Lagrange elements, respectively.}
\label{fig:pentagon_err}
\end{figure}

\begin{example}
\label{ex:fem}
Suppose $\Theta$ is triangulated by $\triangle$ and the Powell--Sabin refinement $\PST$ as in Example~\ref{ex:pentagon}. Moreover, let $\Omega = \Theta$ and $F: \Theta \rightarrow \Omega$ be the identity mapping. We consider the problem \eqref{eq:poisson_problem} with the manufactured solution
\begin{equation}
\label{eq:fem_function}
u(x,y) = \sin(2 \pi x) \cos(2 \pi y),
\end{equation}
shown in Figure~\ref{fig:pentagon_fun}. To compare the Ritz--Galerkin methods for solving this problem based on $\C{0}$ cubic elements and $\C{1}$ cubic Powell--Sabin splines, we generate a series of triangulations starting with $\triangle$ (level 0). Then, for each level, the triangulation from the previous level is refined by the dyadic split, \ie{}, each triangle is divided into four triangles obtained by connecting the midpoints of its edges. Table~\ref{tab:fem} and Figure~\ref{fig:fem_err} provide numerically computed $L^2$ and $L^\infty$ errors and estimated order of convergence with respect to the mesh size $h$. The mesh size is considered as the length of the longest edge of $\triangle$, which is halved on each refinement level.
\end{example}

\begin{table}[h] 
\centering
\begin{tabular}{c|cc|cc|cc|cc|}
 & \multicolumn{4}{c|}{$L^2$ error} & \multicolumn{4}{c|}{$L^\infty$ error} \\ \hline
level & \multicolumn{2}{c|}{Lagrange} & \multicolumn{2}{c|}{Powell--Sabin} & \multicolumn{2}{c|}{Lagrange} & \multicolumn{2}{c|}{Powell--Sabin} \\ \hline
$0$ & $4.84 \cdot 10^{-1}$ & $$ & $5.09 \cdot 10^{-1}$ & $$ & $8.55 \cdot 10^{-1}$ & $$ & $1.02 \cdot 10^{-0}$ & $$ \\
$1$ & $5.79 \cdot 10^{-2}$ & $3.1$ & $3.06 \cdot 10^{-2}$ & $4.1$ & $1.07 \cdot 10^{-1}$ & $3.0$ & $7.94 \cdot 10^{-2}$ & $3.7$ \\
$2$ & $3.96 \cdot 10^{-3}$ & $3.9$ & $2.09 \cdot 10^{-3}$ & $3.9$ & $9.48 \cdot 10^{-3}$ & $3.5$ & $6.92 \cdot 10^{-3}$ & $3.5$ \\
$3$ & $2.52 \cdot 10^{-4}$ & $4.0$ & $1.31 \cdot 10^{-4}$ & $4.0$ & $6.23 \cdot 10^{-4}$ & $3.9$ & $3.84 \cdot 10^{-4}$ & $4.2$ \\
$4$ & $1.57 \cdot 10^{-5}$ & $4.0$ & $8.24 \cdot 10^{-6}$ & $4.0$ & $3.90 \cdot 10^{-5}$ & $4.0$ & $2.35 \cdot 10^{-5}$ & $4.0$
\end{tabular}
\caption{Approximation errors for the problem described in Example~\ref{ex:fem}. In each column section, the value on the left represents the numerically computed error, and the value on the right is the estimated order of convergence with respect to the mesh size.}
\label{tab:fem}
\end{table}

Next, let us consider the cubic splines in the context of isoparametric finite elements. The obtained results show a similar relation between the Powell--Sabin and Lagrange elements as observed in Example~\ref{ex:fem}.

\begin{example}
\label{ex:isofem}
Let $\Theta$ and $\Omega$ be as in Example~\ref{ex:pentagon}, where we consider $\Omega$ parametrized either by the $\C{0}$ geometry mapping $G$ in relation to the cubic Lagrange elements or the $\C{1}$ geometry mapping $F$ in relation to the cubic Powell--Sabin elements. Again, we approximate the manufactured solution \eqref{eq:fem_function} of the problem \eqref{eq:poisson_problem} on a series of triangulations generated as in Example~\ref{ex:fem}. The results are presented in Table~\ref{tab:isofem} and Figure~\ref{fig:isofem_err}.
\end{example}

In the following two examples we test rational Powell--Sabin B-splines on a (curved) planar domain. As in the previous examples, the numerical evidence indicate optimal approximation order.

\begin{figure}[h!]
\centering

\begin{subfigure}[b]{0.48\textwidth}
\centering
\includegraphics[width=\textwidth]{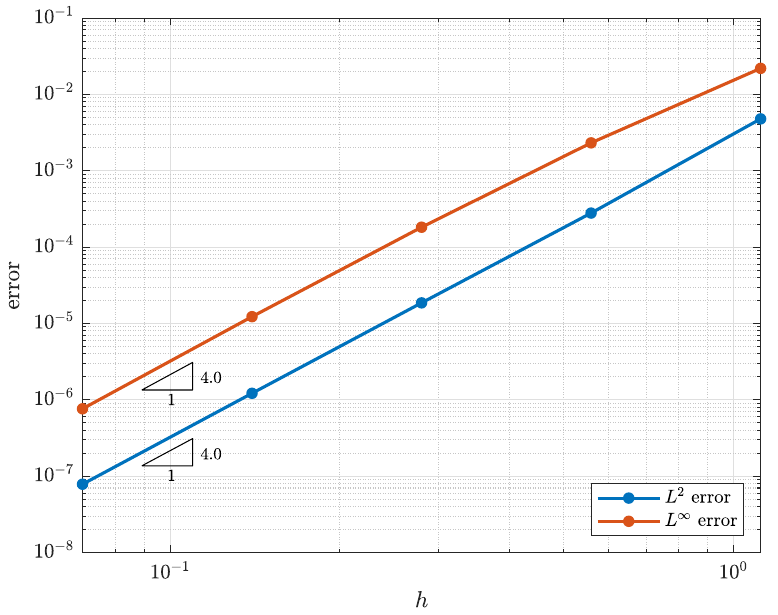}
\caption{Example \ref{ex:poi2}.}
\label{fig:poi2_err}
\end{subfigure}
\hfill
\begin{subfigure}[b]{0.48\textwidth}
\centering
\includegraphics[width=\textwidth]{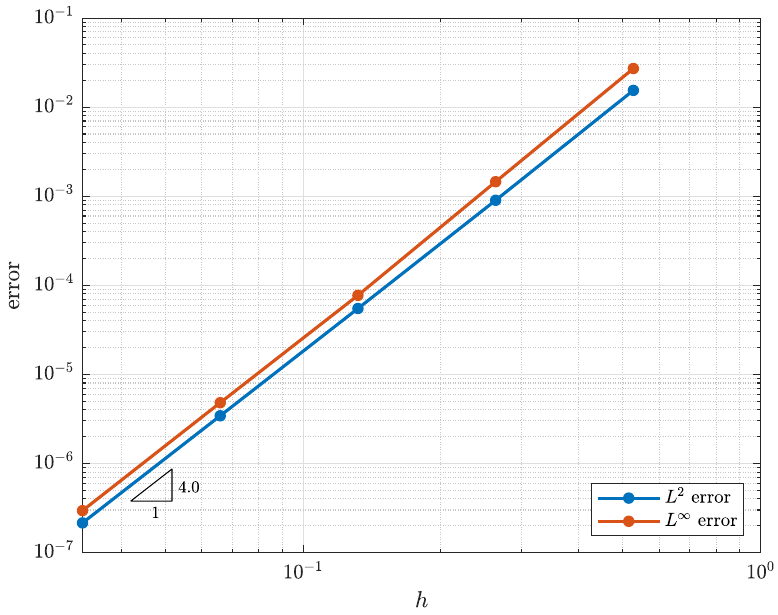}
\caption{Example \ref{ex:poi3}.}
\label{fig:poi3_err}
\end{subfigure}

\caption{Approximation errors for the problems described in Example~\ref{ex:poi2} (left) and Example~\ref{ex:poi3} (right). The errors are presented with respect to the mesh size $h$.}
\label{fig:poi}
\end{figure}

\begin{example}
\label{ex:poi2} 
Suppose the parametric domain $\Theta$ is triangulated by $\triangle$ and the domain $\Omega$ is parametrized by $F: \Theta \rightarrow \Omega$ as in Example~\ref{ex:annulus}. Let us consider the problem \eqref{eq:poisson_problem} with the manufactured solution
\begin{equation}
\label{eq:iga_function}
u(x,y) = \bracket{x^2+\tfrac{1}{3}}e^{2y},
\end{equation}
shown in Figure~\ref{fig:poi2_fun}. We compute approximations from $\mathbb{S}_3^w(\PST)$ on a series of dyadic refinements of $\triangle$. The results are shown in Table~\ref{tab:iga_poi} (left) and Figure~\ref{fig:poi2_err}.
\end{example}

\begin{table}[ht]
\centering
\begin{tabular}{c|cc|cc|cc|cc|}
 & \multicolumn{4}{c|}{$L^2$ error} & \multicolumn{4}{c|}{$L^\infty$ error} \\ \hline
level & \multicolumn{2}{c|}{Lagrange} & \multicolumn{2}{c|}{Powell--Sabin} & \multicolumn{2}{c|}{Lagrange} & \multicolumn{2}{c|}{Powell--Sabin} \\ \hline
$0$ & $6.20 \cdot 10^{-1}$ & $$ & $5.76 \cdot 10^{-1}$ & $$ & $9.56 \cdot 10^{-1}$ & $$ & $9.66 \cdot 10^{-1}$ & $$ \\
$1$ & $8.37 \cdot 10^{-2}$ & $2.9$ & $4.71 \cdot 10^{-2}$ & $3.6$ & $1.67 \cdot 10^{-1}$ & $2.5$ & $1.04 \cdot 10^{-1}$ & $3.2$ \\
$2$ & $6.43 \cdot 10^{-3}$ & $3.7$ & $3.34 \cdot 10^{-3}$ & $3.8$ & $1.47 \cdot 10^{-2}$ & $3.5$ & $1.01 \cdot 10^{-2}$ & $3.4$ \\
$3$ & $4.28 \cdot 10^{-4}$ & $3.9$ & $2.27 \cdot 10^{-4}$ & $3.9$ & $1.21 \cdot 10^{-3}$ & $3.6$ & $7.95 \cdot 10^{-4}$ & $3.8$ \\
$4$ & $2.67 \cdot 10^{-5}$ & $4.0$ & $1.56 \cdot 10^{-5}$ & $3.9$ & $7.96 \cdot 10^{-5}$ & $3.9$ & $5.49 \cdot 10^{-5}$ & $3.9$
\end{tabular}
\caption{Approximation errors for the problem described in Example~\ref{ex:isofem}. In each column section, the value on the left represents the numerically computed error, and the value on the right is the estimated order of convergence with respect to the mesh size.}
\label{tab:isofem}
\end{table}

\begin{table}[h]
\centering
\begin{tabular}{c|cc|cc|}
level & \multicolumn{2}{c|}{$L^2$ error} & \multicolumn{2}{c|}{$L^\infty$ error} \\ \hline
$0$ & $4.79 \cdot 10^{-3}$ & $$ & $2.20 \cdot 10^{-2}$ & $$ \\
$1$ & $2.79 \cdot 10^{-4}$ & $4.1$ & $2.32 \cdot 10^{-3}$ & $3.2$ \\
$2$ & $1.86 \cdot 10^{-5}$ & $3.9$ & $1.82 \cdot 10^{-4}$ & $3.7$ \\
$3$ & $1.22 \cdot 10^{-6}$ & $3.9$ & $1.23 \cdot 10^{-5}$ & $3.9$ \\
$4$ & $7.84 \cdot 10^{-8}$ & $4.0$ & $7.62 \cdot 10^{-7}$ & $4.0$
\end{tabular}
\qquad
\begin{tabular}{c|cc|cc|}
level & \multicolumn{2}{c|}{$L^2$ error} & \multicolumn{2}{c|}{$L^\infty$ error} \\ \hline
$0$ & $1.54 \cdot 10^{-2}$ & $$ & $2.72 \cdot 10^{-2}$ & $$ \\
$1$ & $9.01 \cdot 10^{-4}$ & $4.1$ & $1.46 \cdot 10^{-3}$ & $4.2$ \\
$2$ & $5.49 \cdot 10^{-5}$ & $4.0$ & $7.72 \cdot 10^{-5}$ & $4.2$ \\
$3$ & $3.43 \cdot 10^{-6}$ & $4.0$ & $4.81 \cdot 10^{-6}$ & $4.0$ \\
$4$ & $2.15 \cdot 10^{-7}$ & $4.0$ & $2.95 \cdot 10^{-7}$ & $4.0$
\end{tabular}
\caption{Approximaton errors for the problems described in Example~\ref{ex:poi2} (left) and Example~\ref{ex:poi3} (right). In each column section, the value on the left represents the numerically computed error, and the value on the right is the estimated order of convergence with respect to the mesh size.}
\label{tab:iga_poi}
\end{table}

\begin{figure}[h]
\centering

\begin{subfigure}[b]{0.52\textwidth}
\centering
\begin{subfigure}[b]{0.32\textwidth}
\centering
\includegraphics[width=\textwidth]{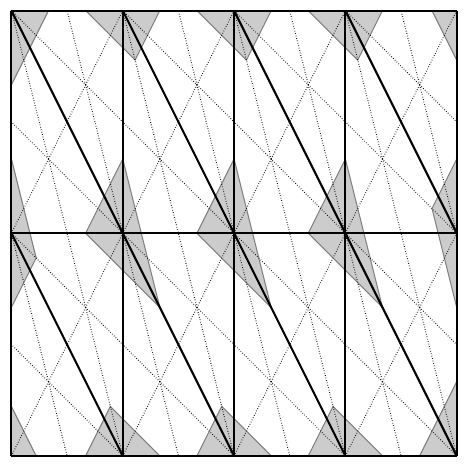}
\caption{Level 1.}
\end{subfigure}
\hfill
\begin{subfigure}[b]{0.32\textwidth}
\centering
\includegraphics[width=\textwidth]{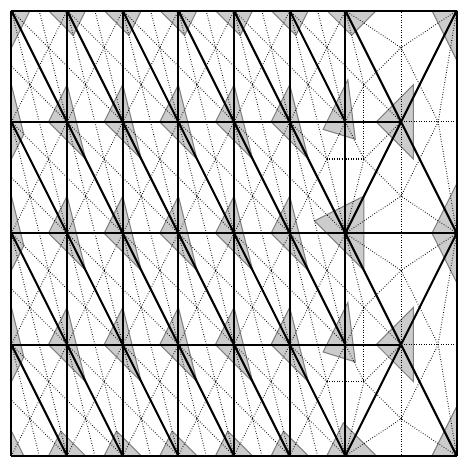}
\caption{Level 2.}
\end{subfigure}
\hfill
\begin{subfigure}[b]{0.32\textwidth}
\centering
\includegraphics[width=\textwidth]{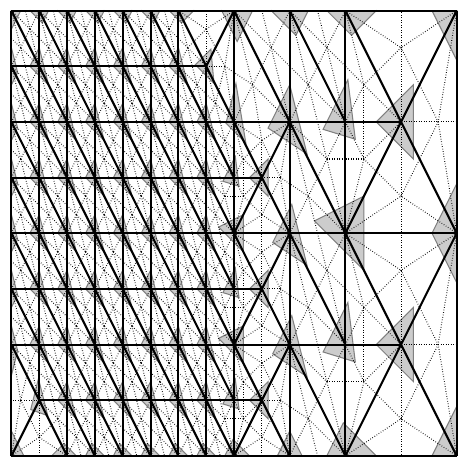}
\caption{Level 3.}
\end{subfigure}

\begin{subfigure}[b]{0.32\textwidth}
\centering
\includegraphics[width=\textwidth]{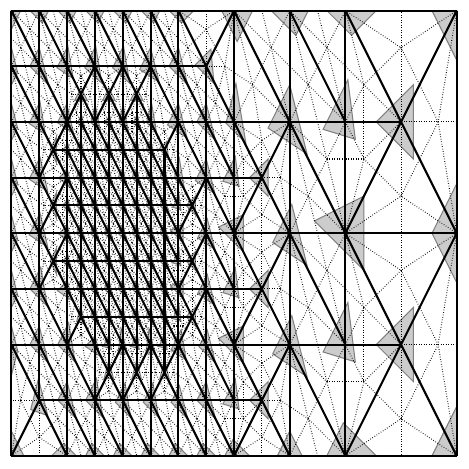}
\caption{Level 4.}
\end{subfigure}
\hfill
\begin{subfigure}[b]{0.32\textwidth}
\centering
\includegraphics[width=\textwidth]{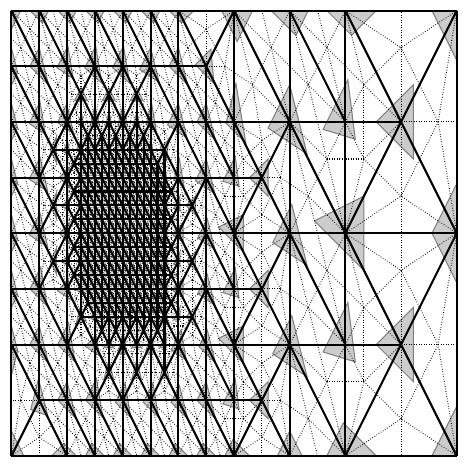}
\caption{Level 5.}
\end{subfigure}
\hfill
\begin{subfigure}[b]{0.32\textwidth}
\centering
\includegraphics[width=\textwidth]{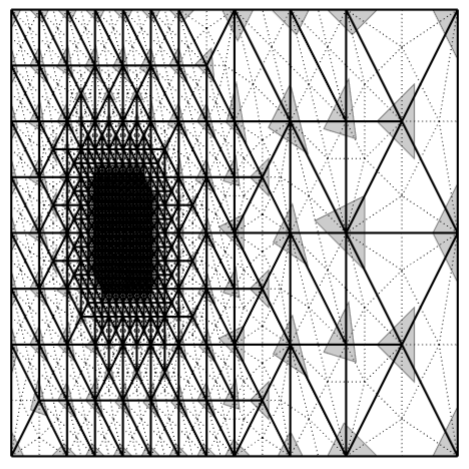}
\caption{Level 6.}
\end{subfigure}
\end{subfigure}
\hfill
\begin{subfigure}[b]{0.47\textwidth}
\centering
\includegraphics[width=\textwidth]{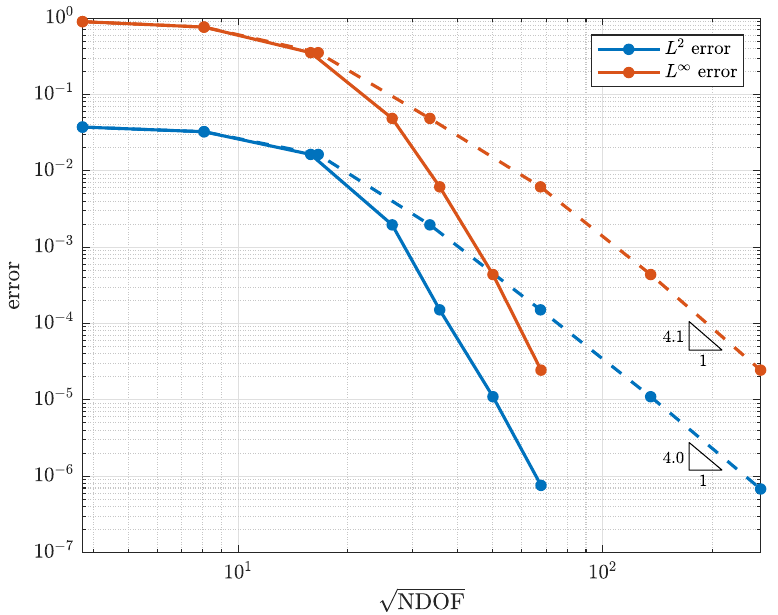}
\caption{Approximation errors.}
\label{fig:poi2_locref_err}
\end{subfigure}

\caption{A series of local refinements of the unit square (left) and approximation errors (right) for the problem described in Example~\ref{ex:poi2_locref}. The errors are presented with respect to $\sqrt{\mathrm{NDOF}}$. The solid and dashed lines correspond to the error obtained on the series of locally and globally refined triangulations, respectively.}
\label{fig:poi2_locref_tri}
\end{figure}

\begin{example}
\label{ex:poi2_locref}
Suppose the setting is the same as in Example~\ref{ex:poi2}, except that this time we consider the problem \eqref{eq:poisson_problem} with the manufactured solution given in polar coordinates $(\rho, \theta)$ by
\begin{equation*}
u(\rho, \theta) = e^{-1000 \bracket{\bracket{\rho-\frac{3}{4}}^2+\bracket{\rho \sin(\theta) - \rho \sin( \frac{1}{8} \pi)}^2}}.
\end{equation*}
Here, $\rho$ is the distance from the origin and $\theta$ is the angle from the polar axis. The graph of $u$ is shown in Figure~\ref{fig:poi2_locref_fun}. Evidently, the graph has a spike at the polar coordinates $(\frac{3}{4}, \frac{1}{8} \pi)$. To take this into account, in addition to the series of dyadic refinements of $\triangle$, we consider the series shown in Figure~\ref{fig:poi2_locref_tri} that conducts local refinement in the region where the function values are subjected to rapid changes. Table~\ref{tab:poi2_locref} and Figure~\ref{fig:poi2_locref_err} provide the comparison of approximation errors obtained by using Powell--Sabin splines on the two series.
\end{example}

\begin{table}[h]
\centering
\begin{tabular}{c|c|c|cc|cc|cc|cc|}
& \multicolumn{2}{c|}{NDOF} & \multicolumn{4}{c|}{$L^2$ error} & \multicolumn{4}{c|}{$L^\infty$ error} \\ \hline
level & global & local & \multicolumn{2}{c|}{global} & \multicolumn{2}{c|}{local} & \multicolumn{2}{c|}{global} & \multicolumn{2}{c|}{local} \\ \hline
$0$ & $14$ & $14$ & $3.74 \cdot 10^{-2}$ & $$ & $3.74 \cdot 10^{-2}$ & $$ & $8.99 \cdot 10^{-1}$ & $$ & $8.99 \cdot 10^{-1}$ & $$ \\
$1$ & $65$ & $65$ & $3.25 \cdot 10^{-2}$ & $0.2$ & $3.25 \cdot 10^{-2}$ & $0.2$ & $7.63 \cdot 10^{-1}$ & $0.2$ & $7.63 \cdot 10^{-1}$ & $0.2$ \\
$2$ & $275$ & $250$ & $1.64 \cdot 10^{-2}$ & $0.9$ & $1.64 \cdot 10^{-2}$ & $1.0$ & $3.53 \cdot 10^{-1}$ & $1.1$ & $3.53 \cdot 10^{-1}$ & $1.1$ \\
$3$ & $1127$ & $700$ & $1.96 \cdot 10^{-3}$ & $3.0$ & $1.96 \cdot 10^{-3}$ & $4.1$ & $4.86 \cdot 10^{-2}$ & $2.8$ & $4.86 \cdot 10^{-2}$ & $3.9$ \\
$4$ & $4559$ & $1276$ & $1.51 \cdot 10^{-4}$ & $3.7$ & $1.51 \cdot 10^{-4}$ & $8.5$ & $6.16 \cdot 10^{-3}$ & $3.0$ & $6.16 \cdot 10^{-3}$ & $6.9$ \\
$5$ & $18335$ & $2500$ & $1.10 \cdot 10^{-5}$ & $3.8$ & $1.10 \cdot 10^{-5}$ & $7.8$ & $4.38 \cdot 10^{-4}$ & $3.8$ & $4.38 \cdot 10^{-4}$ & $7.9$ \\
$6$ & $73535$ & $4588$ & $6.79 \cdot 10^{-7}$ & $4.0$ & $7.57 \cdot 10^{-7}$ & $8.8$ & $2.45 \cdot 10^{-5}$ & $4.2$ & $2.45 \cdot 10^{-5}$ & $9.5$
\end{tabular}
\caption{Approximation errors for the problem described in Example~\ref{ex:poi2_locref}. In each column section related to the errors, the value on the left represents the numerically computed error, and the value on the right is the estimated order of convergence with respect to $\sqrt{\mathrm{NDOF}}$.}
\label{tab:poi2_locref}
\end{table}

Finally, we provide an example of solving the problem \eqref{eq:poisson_problem} on a spatial surface domain. The setup is from \cite[Section 6.1]{iga_dede_15}, where biquadratic and bicubic NURBS are employed to solve the problem. The results reported in Example~\ref{ex:poi3} for Powell--Sabin splines are in line with those obtained by using bicubic NURBS.

\begin{example}
\label{ex:poi3}
Let the parametric domain $\Theta$ be triangulated by $\triangle$ and $\Omega$ parametrized by $F: \Theta \rightarrow \Omega$ as described in Example~\ref{ex:cylinder}. We consider the problem \eqref{eq:poisson_problem} with the manufactured solution given in cylindrical coordinates $(\theta, z)$ by
\begin{equation*}
u(\theta,z) = \bracket{6+4\sqrt{2}} \bracket{1-\cos(\theta)} \bracket{1-\sin(\theta)} \sin\bracket{\tfrac{3}{4} \pi z}.
\end{equation*}
Here, $\theta$ is the azimuth, $z$ is the height, and the radius is omitted as it is constantly equal to $1$. The graph of the function is shown in Figure~\ref{fig:poi3_fun}. On the initial triangulation $\triangle$ we perform a series of dyadic refinements and compute the solutions in the series of spaces $\mathbb{S}_3^w(\PST)$. The results are shown in Table~\ref{tab:iga_poi} (right) and Figure~\ref{fig:poi3_err}.
\end{example}

\subsection{Biharmonic problems}

Let $\Delta_\Omega^2 = \Delta_\Omega \Delta_\Omega$ be the biharmonic operator on the surface domain $\Omega$ parametrized by the geometry mapping $F \in \mathbb{S}_3^w(\PST)^3$. We are interested in finding an approximation $s \in \mathbb{S}_3^w(\PST)$ to a function $u: \Omega \rightarrow \R$ satisfying
\begin{equation}
\label{eq:bih_problem}
\left\{
\begin{aligned}
\Delta_\Omega^2 u &= f && \text{on} \ \Omega \setminus \partial \Omega, \\
u &= g_0 && \text{on} \ \partial \Omega, \\
\nabla_\Omega u \cdot n_{\partial \Omega} &= g_1 && \text{on} \ \partial \Omega,
\end{aligned}
\right.
\end{equation}
for given functions $f: \Omega \rightarrow \R$ and $g_0, g_1: \partial \Omega \rightarrow \R$.  To this end, we follow \cite{iga_bartezzaghi_15} and consider the weak formulation of \eqref{eq:bih_problem} under suitable smoothness conditions on $f$ and $g_0$, $g_1$, in particular $f \in L^2(\Omega)$. After the discretization based on $\mathbb{S}_3^w(\PST)$ and the pull-back to the parametric domain $\Theta$, we arrive at the conditions
\begin{equation}
\label{eq:isofem_conditions_bih}
 \int_\Theta \nabla \cdot \left(\kappa_F \, {K_F}^{-1} \, \grad{s} \right) \, \nabla \cdot \left(\kappa_F \, {K_F}^{-1} \, \grad{N_j}\right) \, \frac{1}{\kappa_F} \, \mathrm{d}\Theta = \int_\Theta (f \circ F) \, N_j \,  \kappa_F \, \mathrm{d}\Theta, \qquad j = 1, 2, \ldots, n_1,
\end{equation}
and the boundary conditions \eqref{eq:boundary_conditions} imposed on the spline $s$ expressed in the form \eqref{eq:spline_approx_boundary} for $d = 1$. Then, by \eqref{eq:isofem_conditions_bih}, the coefficients of $s$ collected in the vector $c = (c_k)_k \in \R^{n_1}$ are the solution of the system of linear equations $A c = b$, where the entries of the matrix $A = (a_{j,k})_{j,k} \in \R^{n_1 \times n_1}$ are given by
\begin{equation*}
a_{j,k} = \int_\Theta \nabla \cdot \left(\kappa_F \, {K_F}^{-1} \, \grad{N_k} \right) \, \nabla \cdot \left(\kappa_F \, {K_F}^{-1} \, \grad{N_j}\right) \, \frac{1}{\kappa_F} \, \mathrm{d}\Theta,
\end{equation*}
and the entries of the vector $b = (b_j)_j \in \R^{n_1}$ are given by
\begin{equation*}
b_j = \int_\Theta \bracket{(f \circ F) \, N_j \, \kappa_F \, - \nabla \cdot \left(\kappa_F \, {K_F}^{-1} \,\grad{s_1}\right) \, \nabla \cdot \left(\kappa_F \, {K_F}^{-1} \, \grad{N_j}\right) \frac{1}{\kappa_F}} \, \mathrm{d}\Theta.
\end{equation*}

In the next examples we demonstrate that the space $\mathbb{S}_3^w(\PST)$ is suitable for computing an approximate solution of the considered fourth order problem.

\begin{figure}[h] 
\centering

\begin{subfigure}[b]{0.48\textwidth}
\centering
\includegraphics[width=\textwidth]{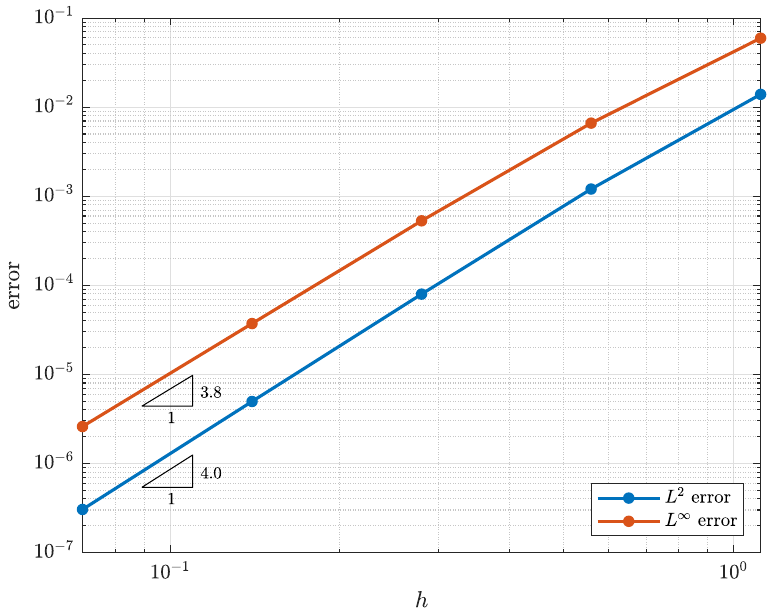}
\caption{Example \ref{ex:bih2}.}
\label{fig:bih2_err}
\end{subfigure}
\hfill
\begin{subfigure}[b]{0.48\textwidth}
\centering
\includegraphics[width=\textwidth]{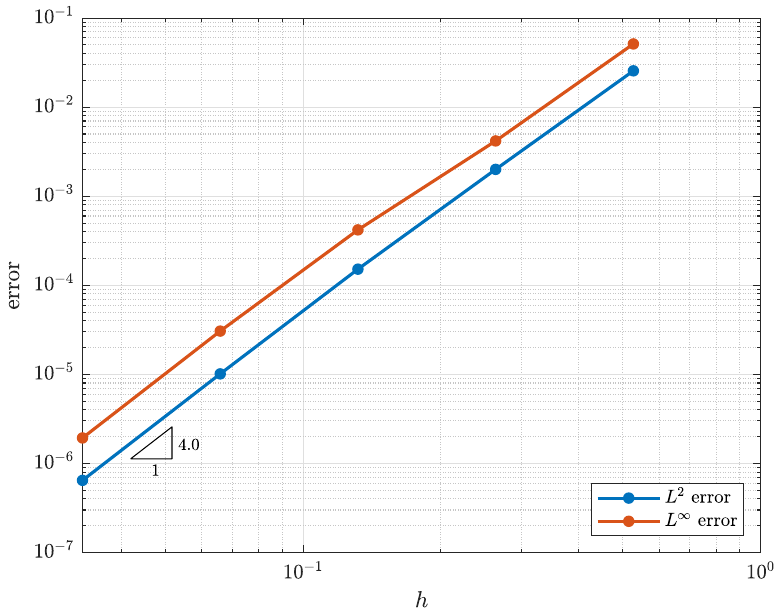}
\caption{Example \ref{ex:bih3}.}
\label{fig:bih3_err}
\end{subfigure}

\caption{Approximation errors for the problems described in Example~\ref{ex:bih2} (left) and Example~\ref{ex:bih3} (right). The errors are presented with respect to the mesh size $h$.}
\label{fig:bih_err}
\end{figure}

\begin{example}
\label{ex:bih2}
Let us repeat Example~\ref{ex:poi2}, this time considering the function $u$ defined by \eqref{eq:iga_function} as the exact solution of the problem \eqref{eq:bih_problem}. The results are shown in Table~\ref{tab:bih} (left) and Figure~\ref{fig:bih2_err}.
\end{example}

\begin{table}[h]
\centering
\begin{tabular}{c|cc|cc|}
level & \multicolumn{2}{c|}{$L^2$ error} & \multicolumn{2}{c|}{$L^\infty$ error} \\ \hline
$0$ & $1.39 \cdot 10^{-2}$ & $$ & $5.96 \cdot 10^{-2}$ & $$ \\
$1$ & $1.20 \cdot 10^{-3}$ & $3.5$ & $6.62 \cdot 10^{-3}$ & $3.2$ \\
$2$ & $7.96 \cdot 10^{-5}$ & $3.9$ & $5.30 \cdot 10^{-4}$ & $3.6$ \\
$3$ & $4.96 \cdot 10^{-6}$ & $4.0$ & $3.72 \cdot 10^{-5}$ & $3.8$ \\
$4$ & $3.04 \cdot 10^{-7}$ & $4.0$ & $2.60 \cdot 10^{-6}$ & $3.8$
\end{tabular}
\qquad
\begin{tabular}{c|cc|cc|}
level & \multicolumn{2}{c|}{$L^2$ error} & \multicolumn{2}{c|}{$L^\infty$ error} \\ \hline
$0$ & $2.56 \cdot 10^{-2}$ & $$ & $5.13 \cdot 10^{-2}$ & $$ \\
$1$ & $2.00 \cdot 10^{-3}$ & $3.7$ & $4.18 \cdot 10^{-3}$ & $3.6$ \\
$2$ & $1.52 \cdot 10^{-4}$ & $3.7$ & $4.19 \cdot 10^{-4}$ & $3.3$ \\
$3$ & $1.01 \cdot 10^{-5}$ & $3.9$ & $3.06 \cdot 10^{-5}$ & $3.8$ \\
$4$ & $6.46 \cdot 10^{-7}$ & $4.0$ & $1.93 \cdot 10^{-6}$ & $4.0$
\end{tabular}
\caption{Approximation errors for the problems described in Example~\ref{ex:bih2} (left) and Example~\ref{ex:bih3} (right). In each column section, the value on the left represents the numerically computed error, and the value on the right is the estimated order of convergence with respect to the mesh size.}
\label{tab:bih}
\end{table}

To conclude, we consider the example conducted in \cite[Section 5.1]{iga_bartezzaghi_15}. This provides another comparison between the Powell--Sabin splines and the biquadratic and bicubic NURBS.

\begin{example}
\label{ex:bih3}
Suppose that the parametric domain $\Theta$ is triangulated by $\triangle$ and $\Omega$ is parametrized by $F: \Theta \rightarrow \Omega$ as described in Example~\ref{ex:cylinder}, where the height of the cylinder shell $h$ is equal to $1$. We aim at finding an approximation to the problem \eqref{eq:bih_problem} with the manufactured solution $u$ given in cylindrical coordinates $(\theta, z)$ by
\begin{equation*}
u(\theta, z) = \sin(2 \theta)^2 \sin(\pi z)^2,
\end{equation*}
shown in Figure~\ref{fig:bih3_fun}. The approximation method is performed on a series of dyadic refinements of $\triangle$. The results are provided in Table~\ref{tab:bih} (right) and Figure~\ref{fig:bih3_err}.
\end{example}

\section{Conclusion}
\label{sec:conclusion}
In this paper, we utilized $\C{1}$ cubic Powell--Sabin splines in the isogeometric analysis framework for the numerical solution of various boundary value problems on planar and spatial surface domains. First, the polynomial $\C{1}$ cubic Powell--Sabin spline space with additional $\C{2}$ super-smoothness properties from \cite{ps3_speleers_15} was reviewed and its rational extension was introduced. Then, both polynomial and rational splines were exploited to define different types of geometry mappings, providing a flexible tool for precise description of (classes of) planar and spatial surfaces. Finally, the obtained surfaces were considered as physical domains in several examples of Poisson and biharmonic problems.

The numerical experiments conducted in this research indicate that $\C{1}$ cubic Powell--Sabin splines perform comparably to or better than cubic Lagrange elements and bicubic NURBS. Their advantage in comparison to Lagrange elements is that they provide $\C{1}$ and not just $\C{0}$ parametrizations. In contrast to tensor product splines, they are not limited to four-sided domains and enable simpler and more versatile domain refinement, which is crucial for efficient approximation. Future work employing higher-degree Powell--Sabin splines \cite{ps_speleers_13,ps_groselj_26} may further enrich the considered approximation framework.

\section*{Acknowledgements}

The work of J.~Gro\v{s}elj and A.~\v{S}adl Praprotnik was partially supported by the ARIS research programme P1-0294. During the preparation of the paper, A.~\v{S}adl Praprotnik was a junior researcher funded by Javna agencija za znanstvenoraziskovalno in inovacijsko dejavnost Republike Slovenije (ARIS).
The work of H.~Speleers was partially supported by the MUR Excellence Department Project MatMod@TOV (CUP E83C23000330006) awarded to the Department of Mathematics of the University of Rome Tor Vergata,
by a Project of Relevant National Interest (PRIN) under the National Recovery and Resilience Plan (PNRR) funded by the European Union -- Next Generation EU (CUP E53D23017910001), and
by the Italian Research Center in High Performance Computing, Big Data and Quantum Computing (CUP E83C22003230001).
H.~Speleers is a member of the research group GNCS (Gruppo Nazionale per il Calcolo Scientifico) of INdAM (Istituto Nazionale di Alta Matematica).

\bibliographystyle{unsrtnat}
\bibliography{references}

@Article{ps3_chen_08,
  author  = {Chen, S. K. and Liu, H. W.},
  title   = {A bivariate ${C}^1$ cubic super spline space on {P}owell--{S}abin triangulation},
  journal = {Comput. Math. Appl.},
  year    = {2008},
  volume  = {56},
  pages   = {1395--1401},
}

@Article{cohen_13,
  author  = {Cohen, E. and Lyche, T. and Riesenfeld, R. F.},
  journal = {Math. Comput.},
  title   = {A {B}-spline-like basis for the {P}owell--{S}abin 12-split based on simplex splines},
  year    = {2013},
  pages   = {1667-1707},
  volume  = {82},
}

@Article{ps_dierckx_97,
  Title                    = {On calculating normalized {P}owell--{S}abin {B}-splines},
  Author                   = {Dierckx, P.},
  Journal                  = {Comput. Aided Geom. Design},
  Year                     = {1997},
  Pages                    = {61--78},
  Volume                   = {15}
}

@Article{iga_eddargani_24,
  Title                    = {Quadrature rules for {$C^1$} quadratic spline finite elements on the {P}owell--{S}abin 12-split},
  Author                   = {Eddargani, S. and Lyche, T. and Manni, C. and Speleers, H.},
  Journal                  = {Comput. Methods Appl. Mech. Engrg.},
  Year                     = {2024},
  Pages                    = {117196},
  Volume                   = {430}
}

@Article{farin_86,
  Title                    = {Triangular {B}ernstein--{B}\'{e}zier patches},
  Author                   = {Farin, G.},
  Journal                  = {Comput. Aided Geom. Design},
  Year                     = {1986},
  Pages                    = {83--127},
  Volume                   = {3}
}

@Article{ps_groselj_16,
  author  = {Gro\v{s}elj, J.},
  title   = {A normalized representation of super splines of arbitrary degree on {P}owell--{S}abin triangulations},
  journal = {BIT Numer. Math.},
  year    = {2016},
  volume  = {56},
  pages   = {1257--1280},
}

@Article{ps3_groselj_16,
  author  = {Gro\v{s}elj, J. and Krajnc, M.},
  title   = {${C}^1$ cubic splines on {P}owell--{S}abin triangulations},
  journal = {Appl. Math. Comput.},
  year    = {2016},
  volume  = {272},
  pages   = {114--126},
}

@Book{lai_07,
  title     = {{S}pline {F}unctions on {T}riangulations},
  publisher = {Cambridge University Press},
  year      = {2007},
  author    = {Lai, M.-J. and Schumaker, L. L.},
}

@Article{ps3_lamnii_14,
  author  = {Lamnii, M. and Mraoui, H. and Tijini, A. and Zidna, A.},
  title   = {A normalized basis for {$C^1$} cubic super spline space on {P}owell--{S}abin triangulation},
  journal = {Math. Comput. Simul.},
  year    = {2014},
  volume  = {99},
  pages   = {108--124},
}

@Article{lyche_25,
  author  = {Lyche, T. and Merrien, J.-L. and Speleers, H.},
  journal = {Comput. Aided Geom. Design},
  title   = {A {$C^1$} simplex-spline basis for the {A}lfeld split in {${\mathbb R}^s$}},
  year    = {2025},
  pages   = {102412},
  volume  = {117},
}

@Article{speleers_06,
  Title                    = {Numerical solution of partial differential equations with {P}owell--{S}abin splines},
  Author                   = {Speleers, H. and Dierckx, P. and Vandewalle, S.},
  Journal                  = {J. Comput. Appl. Math.},
  Year                     = {2006},
  Pages                    = {643--659},
  Volume                   = {189}
}

@Article{ps5_speleers_10,
  Title                    = {A normalized basis for quintic {P}owell--{S}abin splines},
  Author                   = {Speleers, H.},
  Journal                  = {Comput. Aided Geom. Design},
  Year                     = {2010},
  Pages                    = {438--457},
  Volume                   = {27}
}

@Article{rct3_speleers_10,
  Title                    = {A normalized basis for reduced {C}lough--{T}ocher splines},
  Author                   = {Speleers, H.},
  Journal                  = {Comput. Aided Geom. Design},
  Year                     = {2010},
  Pages                    = {700--712},
  Volume                   = {27}
}

@Article{speleers_12,
  Title                    = {Isogeometric analysis with {P}owell--{S}abin splines for advection-diffusion-reaction problems},
  Author                   = {Speleers, H. and Manni, C. and Pelosi, F. and Sampoli, M. L.},
  Journal                  = {Comput. Methods Appl. Mech. Engrg.},
  Year                     = {2012},
  Pages                    = {132--148},
  Volume                   = {221--222}
}

@Article{ps_speleers_13,
  Title                    = {Construction of normalized {B}-splines for a family of smooth spline spaces over {P}owell--{S}abin triangulations},
  Author                   = {Speleers, H.},
  Journal                  = {Constr. Approx.},
  Year                     = {2013},
  Pages                    = {41--72},
  Volume                   = {37}
}

@Article{ps3_speleers_15,
  Title                    = {A new {B}-spline representation for cubic splines over {P}owell--{S}abin triangulations},
  Author                   = {Speleers, H.},
  Journal                  = {Comput. Aided Geom. Design},
  Year                     = {2015},
  Pages                    = {42--56},
  Volume                   = {37}
}

@Article{ps3_groselj_17,
  author  = {Gro\v{s}elj, J. and Speleers, H.},
  title   = {Construction and analysis of cubic {P}owell--{S}abin {B}-splines},
  journal = {Comput. Aided Geom. Design},
  year    = {2017},
  volume  = {57},
  pages   = {1--22},
}

@Book{fem_ciarlet_02,
  title     = {The Finite Element Method for Elliptic Problems},
  publisher = {Society for Industrial and Applied Mathematics},
  year      = {2002},
  author    = {Ciarlet, P. G.},
  edition   = {2nd},
}

@Article{iga_hughes_05,
  author  = {Hughes, T. J. R. and Cottrell, J. A. and Bazilevs, Y.},
  title   = {Isogeometric analysis: {CAD}, finite elements, {NURBS}, exact geometry and mesh refinement},
  journal = {Comput. Methods Appl. Mech. Engrg.},
  year    = {2005},
  volume  = {194},
  pages   = {4135--4195},
}

@Article{iga_jaxon_14,
  author  = {Jaxon, N. and Qian, X.},
  title   = {Isogeometric analysis on triangulations},
  journal = {Comput. Aided Design},
  year    = {2014},
  volume  = {46},
  pages   = {45--57},
}

@Article{iga_xia_15,
  author  = {Xia, S. and Wang, X. and Qian, X.},
  title   = {Continuity and convergence in rational triangular {B}\'{e}zier spline based isogeometric analysis},
  journal = {Comput. Methods Appl. Mech. Engrg.},
  year    = {2015},
  volume  = {297},
  pages   = {292--324},
}

@Article{iga_wang_18,
  author  = {Wang, C. and Xia, S. and Wang, X. and Qian, X.},
  title   = {Isogeometric shape optimization on triangulations},
  journal = {Comput. Methods. Appl. Mech. Engrg},
  year    = {2018},
  volume  = {331},
  pages   = {585--622},
}

@Article{iga_xia_18,
  author  = {Xia, S. and Qian, X.},
  title   = {Generating high-quality high-order parametrization for isogeometric analysis on triangulations},
  journal = {Comput. Methods Appl. Mech. Engrg.},
  year    = {2018},
  volume  = {338},
  pages   = {1--26},
}

@Article{iga_zareh_19,
  author  = {Zareh, M. and Qian, X.},
  title   = {{K}irchhoff--{L}ove shell formulation based on triangular isogeometric analysis},
  journal = {Comput. Methods Appl. Mech. Engrg.},
  year    = {2019},
  volume  = {347},
  pages   = {853--873},
}

@Article{iga_giorgiani_18,
  author  = {Giorgiani, G. and Guillard, H. and Nkonga, B. and Serre, E.},
  title   = {A stabilized {P}owell--{S}abin finite-element method for the {2D} {E}uler equations in supersonic regime},
  journal = {Comput. Methods Appl. Mech. Engrg.},
  year    = {2018},
  volume  = {340},
  pages   = {216--235},
}

@Article{bell_groselj_20,
  author  = {Gro\v{s}elj, J. and Knez, M.},
  title   = {On stable representations of {B}ell elements},
  journal = {Comput. Math. Appl.},
  year    = {2020},
  volume  = {79},
  pages   = {2924--2941},
}

@Article{rps_speleers_13,
  author  = {Speleers, H. and Manni, C. and Pelosi, F.},
  title   = {From {NURBS} to {NURPS} geometries},
  journal = {Comput. Methods Appl. Mech. Engrg.},
  year    = {2013},
  volume  = {255},
  pages   = {238--254},
}

@Article{ps2_speleers_15,
  author  = {Speleers, H. and Manni, C.},
  journal = {J. Comput. Appl. Math.},
  title   = {Optimizing domain parameterization in isogeometric analysis based on {P}owell--{S}abin splines},
  year    = {2015},
  pages   = {68--86},
  volume  = {289},
}

@Article{rps_speleers_07,
  author  = {Speleers, H. and Dierckx, P. and Vandewalle, S.},
  title   = {Weight control for modelling with {NURPS} surfaces},
  journal = {Comput. Aided Geom. Design},
  year    = {2007},
  volume  = {24},
  pages   = {179--186},
}

@Article{ps3_groselj_18,
  author  = {Gro\v{s}elj, J. and Speleers, H.},
  title   = {Three recipes for quasi-interpolation with cubic {P}owell--{S}abin splines},
  journal = {Comput. Aided Geom. Design},
  year    = {2018},
  volume  = {67},
  pages   = {47--70},
}

@Article{arg_groselj_21,
  author  = {Gro\v{s}elj, J. and Knez, M.},
  title   = {A construction of edge {B}-spline functions for a {$C^1$} polynomial spline on two triangles and its application to {A}rgyris type splines},
  journal = {Comput. Math. Appl.},
  year    = {2021},
  volume  = {99},
  pages   = {329--344},
}

@Article{thb_giannelli_12,
  author  = {Giannelli, C. and J\"{u}ttler, B. and Speleers, H.},
  title   = {{THB}-splines: The truncated basis for hierarchical splines},
  journal = {Comput. Aided Geom. Design},
  year    = {2012},
  volume  = {29},
  number  = {7},
  pages   = {485--498},
}

@Article{ts_sederberg_03,
  author  = {Sederberg, T. W. and Zheng, J. and Bakenov, A. and Nasri, A.},
  title   = {{T}-splines and {T}-{NURCC}s},
  journal = {ACM Trans. Graph.},
  year    = {2003},
  volume  = {22},
  pages   = {477--484},
}

@Article{lr_dokken_13,
  author  = {Dokken, T. and Lyche, T. and Pettersen, K. F.},
  title   = {Polynomial splines over locally refined box-partitions},
  journal = {Comput. Aided Geom. Design},
  year    = {2013},
  volume  = {30},
  pages   = {331--356},
}

@Book{iga_cottrell_09,
  title     = {Isogeometric Analysis: Toward Integration of {CAD} and {FEA}},
  publisher = {John Wiley \& Sons},
  year      = {2009},
  author    = {Cottrell, J. A. and Hughes, T. J. R. and Bazilevs, Y.},
}

@Article{speleers_07,
  author  = {Speleers, H. and Dierckx, P. and Vandewalle, S.},
  title   = {{P}owell--{S}abin splines with boundary conditions for polygonal and non-polygonal domains},
  journal = {J. Comput. Appl. Math.},
  year    = {2007},
  volume  = {206},
  number  = {1},
  pages   = {55--72},
}

@Article{rps3_groselj_25,
  author  = {Gro\v{s}elj, J. and \v{S}adl Praprotnik, A.},
  journal = {J. Comput. Appl. Math.},
  title   = {Rational {$C^1$} cubic {P}owell--{S}abin {B}-splines with application to representation of ruled surfaces},
  year    = {2025},
  volume  = {457},
}

@Article{iga_dede_15,
  author  = {Ded\`{e}, L. and Quarteroni, A.},
  journal = {Comput. Methods Appl. Mech. Engrg.},
  title   = {Isogeometric analysis for second order partial differential equations on surfaces},
  year    = {2015},
  pages   = {807--834},
  volume  = {284},
}

@Article{iga_bartezzaghi_15,
  author  = {Bartezzaghi, A. and Ded\`{e}, L. and Quarteroni, A.},
  journal = {Comput. Methods. Appl. Mech. Engrg.},
  title   = {Isogeometric analysis of high order partial differential equations on surfaces},
  year    = {2015},
  pages   = {446--469},
  volume  = {295},
}

@Article{kestelman_71,
  author  = {Kestelman, H.},
  journal = {Amer. Math. Monthly},
  title   = {Mappings with non-vanishing {J}acobian},
  year    = {1971},
  pages   = {662--663},
  volume  = {78},
}

@Book{deboor_01,
  author    = {de Boor, C.},
  publisher = {Springer New York},
  title     = {A Practical Guide to Splines},
  year      = {2001},
  series    = {Applied Mathematical Sciences},
}

@Article{ct_groselj_22,
  author  = {Gro\v{s}elj, J. and Knez, M.},
  journal = {Comput. Methods Appl. Mech. Engrg.},
  title   = {Generalized {$C^1$} {C}lough--{T}ocher splines for {CAGD} and {FEM}},
  year    = {2022},
  pages   = {114983},
  volume  = {395},
}

@InProceedings{ps_windmolders_99,
  author    = {Windmolders, J. and Dierckx, P.},
  booktitle = {Proceedings of Curve and Surface Fitting},
  title     = {From {PS}-splines to {NURPS}},
  year      = {1999},
  editor    = {Cohen, A. and Rabut, C.},
  pages     = {45--54},
  publisher = {Vanderbilt University Press},
}

@InProceedings{ps_windmolders_00,
  author    = {Windmolders, J. and Dierckx, P.},
  booktitle = {Proceedings of Mathematical Methods for Curves and Surfaces},
  title     = {{NURPS} for special effects and quadrics},
  year      = {2000},
  editor    = {Lyche, T. and Schumaker, L. L.},
  pages     = {527--534},
  publisher = {Vanderbilt University Press},
}

@InProceedings{liu_07,
  author    = {Liu, Y. and Snoeyink, J.},
  booktitle = {Proc. of the 23rd Annual Symposium on Computational Geometry},
  title     = {Quadratic and cubic {B}-splines by generalizing higher-order {V}oronoi diagrams},
  year      = {2007},
  pages     = {150--157},
  publisher = {ACM},
}

@Article{lyche_18,
  author  = {Lyche, T. and Merrien, J.-L.},
  journal = {Comput. Aided Geom. Design},
  title   = {Simplex-splines on the {C}lough--{T}ocher element},
  year    = {2018},
  pages   = {76--92},
  volume  = {65},
}

@Article{lyche_22,
  author  = {Lyche, T. and Manni, C. and Speleers, H.},
  journal = {Found. Comput. Math.},
  title   = {Construction of {$C^2$} cubic splines on arbitrary triangulations},
  year    = {2022},
  pages   = {1309--1350},
  volume  = {22},
}

@Article{wang_22,
  author  = {Wang, Z. and Cao, J. and Wei, X. and Zhang, Y. J.},
  journal = {Comput. Methods Appl. Mech. Engrg.},
  title   = {{TCB}-spline-based isogeometric analysis method with high-quality parameterizations},
  year    = {2022},
  pages   = {114771},
  volume  = {393},
}

@Article{lyche_24,
  author  = {Lyche, T. and Manni, C. and Speleers, H.},
  journal = {Appl. Math. Comput.},
  title   = {A local simplex spline basis for {$C^3$} quartic splines on arbitrary triangulations},
  year    = {2024},
  volume  = {462},
}

@Book{Delfour2011,
  author    = {Delfour, M. C. and Zol\'{e}sio, J. P.},
  publisher = {Society for Industrial and Applied Mathematics},
  title     = {Shapes and Geometries: Metrics, Analysis, Differential Calculus, and Optimization},
  year      = {2011},
  edition   = {2nd},
  series    = {Advances in Design and Control},
}

@Article{ps_groselj_26,
  author  = {J. Gro\v{s}elj},
  journal = {Math. Comput. Simul.},
  title   = {Higher-degree super-smooth ${C}^1$ splines over a {P}owell--{S}abin refined triangulation},
  year    = {2026},
  pages   = {382--406},
  volume  = {243},
}

\end{document}